\documentclass[preprint]{ourstyle}

\usepackage[american]{babel}
\usepackage{bm}
\usepackage{amsmath}
\usepackage{amssymb}
\usepackage{amsbsy}
\usepackage{amsfonts}
\usepackage{array}
\usepackage{verbatim}
\usepackage{graphicx}
\usepackage{epstopdf,pdfsync}
\usepackage{pdfsync}
\usepackage{listings}
\usepackage{subfigure}
\usepackage{soul}
\usepackage{multicol}

\usepackage{tikz}
\usetikzlibrary{shapes,arrows}
\usepackage[noprefix]{nomencl}
\usepackage{booktabs}
\usepackage{stmaryrd}

\graphicspath{{figure/}}
\DeclareGraphicsExtensions{.pdf,.jpg,.png,.eps}

\usepackage{color}
\usepackage{afterpage}

\pgfdeclarelayer{background}
\pgfdeclarelayer{foreground}
\pgfsetlayers{background,main,foreground}

\renewcommand{\emptyset}{\varnothing}
\newcommand{\dsum}{\displaystyle\sum}

\newcommand{\dinf}{\displaystyle\inf}

\newcommand{\dint}{\displaystyle\int}
\newcommand{\grad}{\nabu}
\newcommand{\Div}{{\rm div}}
\newcommand{\nabu}{{\boldsymbol \nabla}}
\newcommand{\Frac}{\displaystyle\frac}

\renewcommand{\vector}[1]{\boldsymbol{#1}}

\providecommand{\norm}[1]{\left\|{#1}\right\|}
\providecommand{\abs}[1]{\left|{#1}\right|}

\newtheorem{remark}{Remark}

\newtheorem{proposition}{Proposition}

\newcommand{\NEW}[1]{#1}

\begin{document}

\begin{frontmatter}

\title{A Multiscale Thermo-Fluid Computational Model \\
for a Two-Phase Cooling System}

\address[POLI]{Dipartimento di Matematica \lq\lq F.Brioschi\rq\rq,
               Politecnico di Milano, \\
	           Piazza L. da Vinci 32, 20133 Milano Italy}

\address[Indiana]{IUPUI
Department of Mathematical Sciences \\
402 N. Blackford St., LD 270 E
Indianapolis, IN 46202-3267}

\address[ABB]{ABB Switzerland Ltd., Corporate Research\\
              Segelhofstrasse 1K 5405, Baden-D\"attwil, Switzerland}

\address[POLI_MOX]{MOX Modeling and Scientific Computing, \\ 
                   Dipartimento di Matematica \lq\lq F.Brioschi\rq\rq,
                   Politecnico di Milano, \\
                   Piazza L. da Vinci 32, 20133 Milano Italy}

\address[CEN]{CEN Centro Europeo di Nanomedicina,\\
              Piazza L. da Vinci 32, 20133 Milano Italy}


\author[POLI]{Riccardo Sacco}
\author[Indiana]{Lucia Carichino}
\author[POLI_MOX,CEN]{Carlo de Falco}
\author[POLI]{Maurizio Verri}
\author[ABB]{Francesco Agostini}
\author[ABB]{Thomas Gradinger}

\date{\today}
\journal{Comp. Meth. Appl. Mech. Engrg.}

\begin{abstract}

In this paper, we describe a mathematical model and 
a numerical simulation method for the condenser component
of a novel two-phase thermosyphon cooling system
for power electronics applications.
The condenser consists of a set of 
roll-bonded vertically mounted fins among which air flows 
by either natural or forced convection.
In order to deepen the understanding of the mechanisms that determine 
the performance of the condenser and to facilitate the further 
optimization of its industrial design, a multiscale approach 
\NEW{is} developed to reduce as much as possible the complexity of 
the simulation code while maintaining reasonable predictive accuracy.
To this end, heat diffusion in the fins and its convective 
transport in air are modeled as 2D processes while
the flow of the two-phase coolant within the fins 
is modeled as a 1D network of pipes.
For the numerical solution of the resulting equations, 
a Dual Mixed-Finite Volume scheme with 
Exponential Fitting stabilization is used for 2D heat diffusion 
and convection while
a Primal Mixed Finite Element discretization method with upwind stabilization
 is used for the 1D coolant flow. 
The mathematical model and the numerical method are validated 
through extensive simulations of realistic device structures
which prove to be in excellent agreement with available 
experimental data.

\end{abstract}


\end{frontmatter}

\section{Introduction and Motivation}\label{sec:intro}
\begin{figure}[h!]
\centering
\includegraphics[width=0.35\columnwidth]{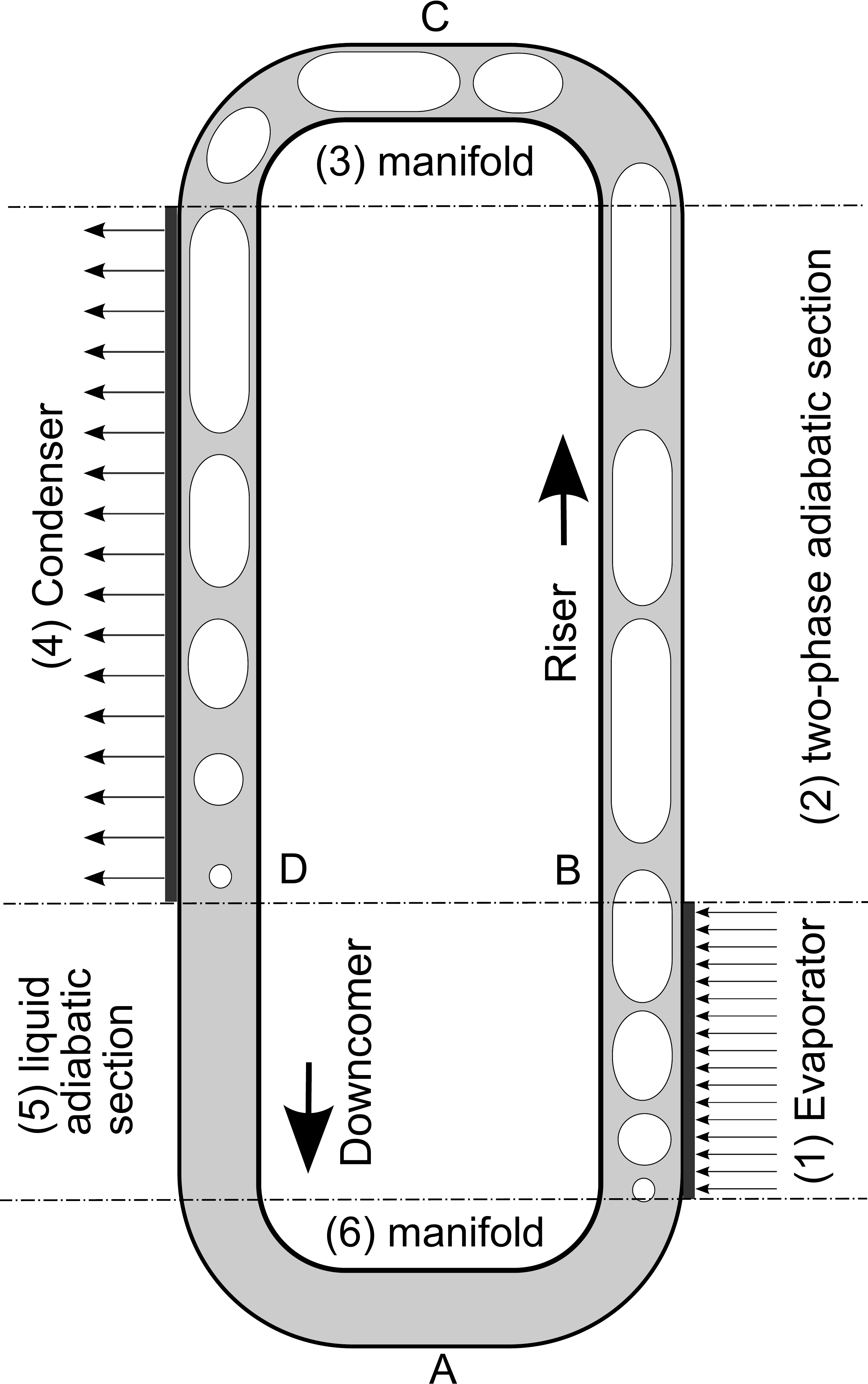}
\caption{Schematic representation of the working principle of a thermosyphon cooler.}
\label{fig:tpcoolermod}
\end{figure}

\begin{figure}[h!]
\begin{minipage}[c]{0.49\linewidth}
\centering
\subfigure[Cooling system assembly]{
\includegraphics[width=0.99\columnwidth]{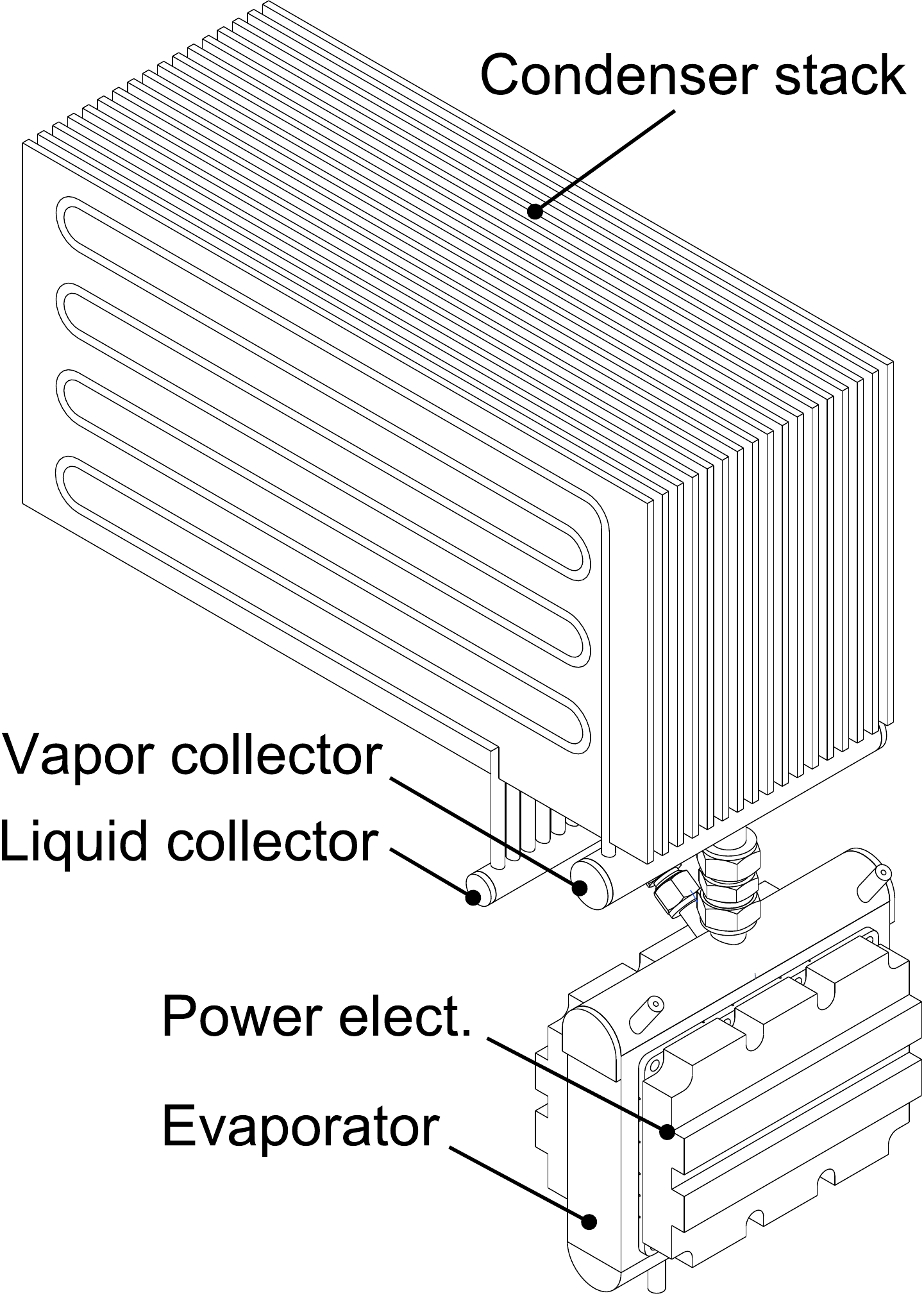}
\label{fig:tpcooler}}
\end{minipage}
\quad
\begin{minipage}[c]{0.49\linewidth}
\subfigure[Detail of the condenser and symmetry plane]{
\includegraphics[width=0.99\columnwidth]{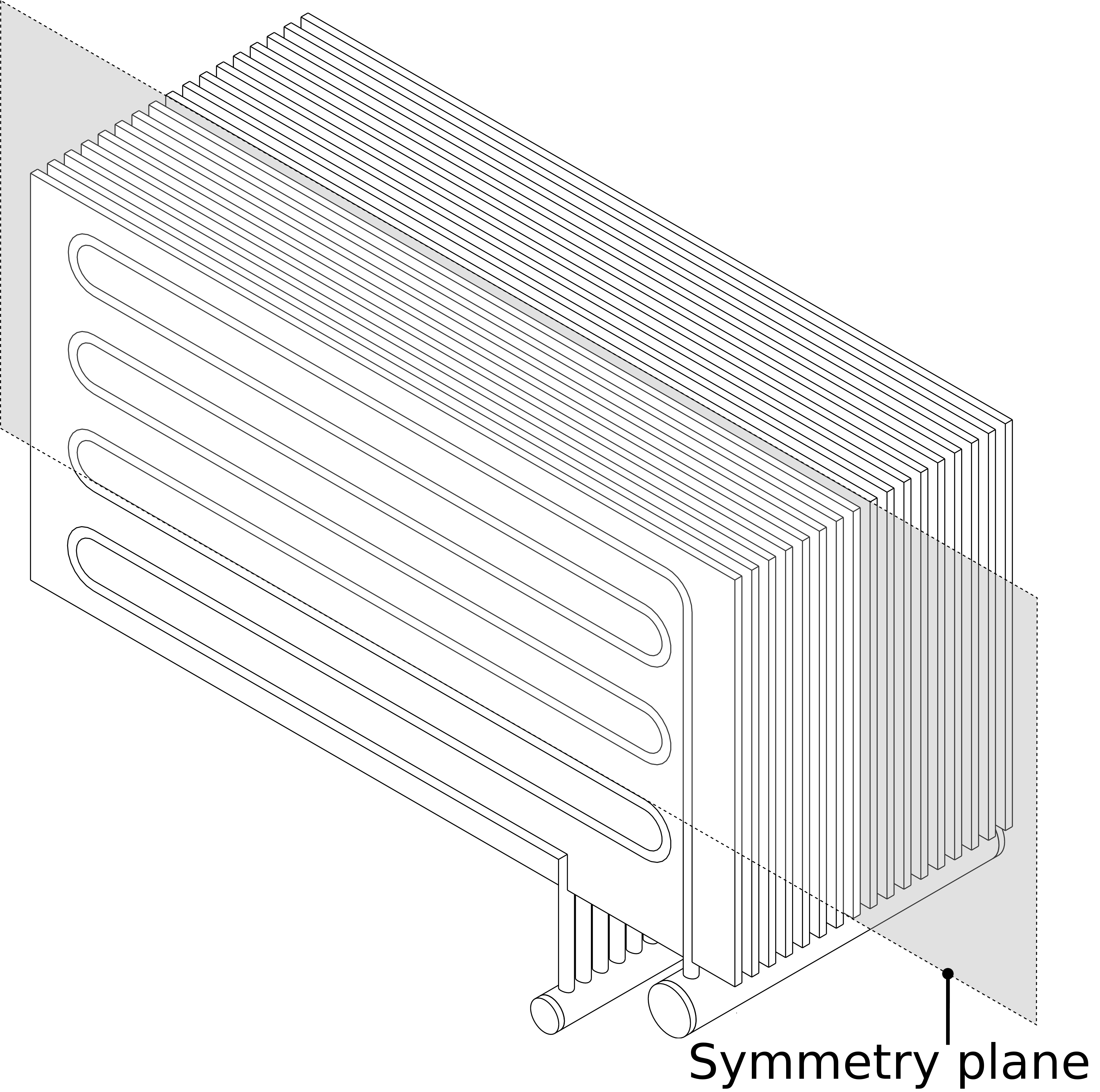}
\label{fig:roll1}
}
\end{minipage}
\caption{Two-phase cooling system based on the thermosyphon principle~\cite{experimental,iecon11}.}
\end{figure}

Ever since the early 1980s the increasing growth of new technologies and applications 
has been shifting scientific interest on power electronics. 
In such wide-range industrial context, the necessity to develop devices with a high power dissipation per unit volume 
has justified the need of advanced cooling systems
capable to prevent  excessive temperature increase
and consequent device failure. Conventional cooling procedures 
exploit convection heat transfer between a fluid in motion and a 
bounding surface at different temperatures. Typical examples 
are water-cooled and air-cooled systems, widely used in power electronics 
applications. A different approach to cooling is represented by 
the two-phase thermosyphon device \NEW{whose functioning principle is schematically illustrated} in Fig.~\ref{fig:tpcoolermod}
\NEW{and whose structure is shown in Fig.~\ref{fig:tpcooler}.}

This kind of device consists of an evaporator and a condenser 
connected by a pipe in which a mixture of liquid and vapor phases is flowing. 
The heat generated by an electronic device
in contact with the evaporator is collected by means of an evaporating fluid. 
The vapor phase fluid, rising in the pipe, passes through the condenser where 
it returns to the liquid phase. As no pumps are needed to move the refrigerant 
fluid from the evaporator to the condenser,
the resulting thermodynamical efficiency of two-phase cooling systems is 
remarkably superior to that of water-cooled or air-cooled systems 
(see~\cite{Lasance1997}).
In order to deepen our understanding of the mechanisms that determine 
the performance of a two-phase thermosyphon cooler device
and to facilitate the further optimization 
of its design, in the present research we focus 
on the study of the condenser subsystem (see Fig.~\ref{fig:roll1}),
for which we \NEW{develop} a multiscale mathematical model that \NEW{is}
implemented in a numerical simulation code.
As computational efficiency is a stringent requirement in
industrial design and optimization procedures, 
model complexity is suitably reduced through the adoption of 
physically sound consistent assumptions that allow \NEW{us} to end up with
a system of nonlinearly coupled 2D PDEs for the air and panel 
temperatures, and 1D PDEs for the refrigerant fluid flow.

Another important constraint is represented by the ability
of the computational method to reproduce on the discrete level 
important physical features characterizing the problem at hand, 
such as mass and flux conservation, and its robustness in the presence of 
dominating convective flow regimes. 
These requirements are here satisfied by the introduction of
a stabilized mixed finite element scheme on quadrilateral 
grids that automatically provides the desired 
inter-element flux conservation and upwinding through the
use of suitable quadrature rules for the mass flux matrix 
and convection term. The resulting discrete method has also
an immediate interpretation in terms of finite volume formulation 
which allows a compact implementation of the scheme
that highly improves the overall efficiency of
the computer-aided design procedure.

A final issue of critical importance in the development of
a reliable computational tool for use in industrial design
is model calibration and validation. 
Model calibration is properly addressed by supplying
the parameter setting in the equation system with 
suitable {\em empirical correlations}, that are functional relations 
between two or more physical variables, usually obtained by means of a series 
of experimental tests. In common engineering practice, 
correlations are widely used because they allow to account for 
complex physical phenomena in a simple and synthetic manner,
albeit their applicatibility is clearly restricted to a specific
admissible range of parameter values. 
Model validation is carried out through extensive numerical simulations 
of the two-phase condenser under realistic working conditions.

An outline of the article is as follows. 
Sect.~\ref{sec:mathmod2d} describes the two--dimensional model 
for heat convection in air and heat diffusion in the panel
whose derivation from the corresponding full 3D model
is outlined in~\ref{sec:3d_2_2d}.
The simplified geometrical representation of the coolant--filled
channel and the one--dimensional system of PDEs describing the 
flow within it are dealt with in Sect.~\ref{sec:mathmod1d}.
Sect.~\ref{sec:iter} discusses the decoupled iterative algorithm used
to solve the complete model \NEW{while Sects.~\ref{sec:air_wall_discr} 
and~\ref{sec:chan_discr}} are devoted to the
discussion of the discretization techniques adopted to 
treat each differential subsystems arising from system linearization.
Finally, in Sect.~\ref{sec:results} simulation
results are \NEW{presented and discussed} and in Sect.~\ref{sec:conclusions} conclusions are
drawn \NEW{and possible future research directions are addressed.}

\section{Mathematical Models}

In this section we describe the mathematical model on which our numerical 
simulation tool for the condenser is based. The equations for heat 
convection in air and heat diffusion in the panel wall are presented in 
Sect.~\ref{sec:mathmod2d}, while the model for the two-phase flow 
in the channel is in Sect.~\ref{sec:channel_domain}. 

\subsection{2D model for the panel wall and air flow}
\label{sec:mathmod2d}

The model for heat diffusion and convection 
is based on the following set of simplifying assumptions:
\begin{enumerate}
\item [(H1)] the geometry of the channel embedded into each panel 
of the condenser is the same;
\item [(H2)] air flow is in steady-state conditions;
\item [(H3)] air flow conditions in between each pair of condenser fins are identical;
\item [(H4)] air flow velocity $\vector{v}_{a}$ is everywhere parallel to the fin plates and its magnitude varies only in the orthogonal direction; 
\item [(H5)] air density $\rho_{a}$ is constant; 
\item [(H6)] the thickness of each panel is negligible compared to its size in any other direction;
\item [(H7)] 
the thickness of the air layer separating two panels 
in the condenser is negligible compared to the panel size in any other direction.
\end{enumerate}
\begin{figure}[h!]
\centering
\includegraphics[width=.75\linewidth]{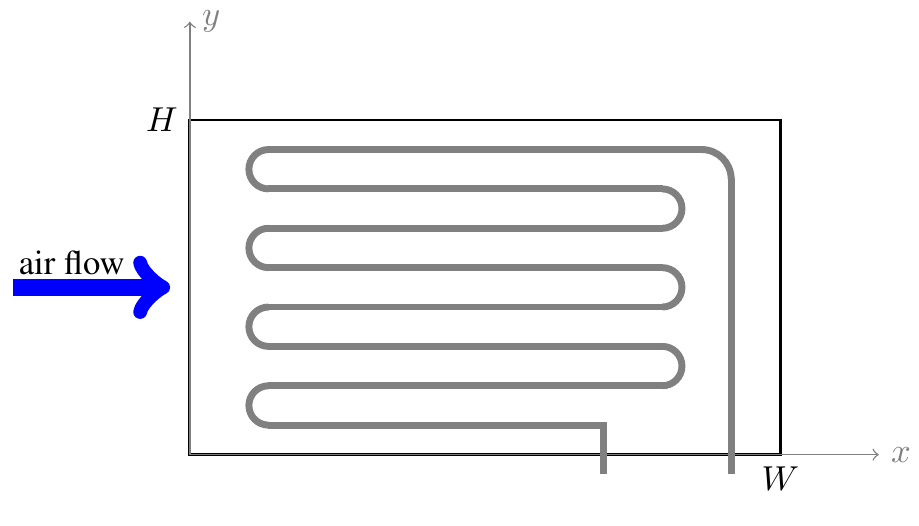}
\caption{Computational domain for the two-dimensional heat flow problem.}
\label{fig:symmetry}
\end{figure}

Under the assumptions above, symmetry considerations lead to 
define the simplified computational domain 
$\Omega:= (0,W)\times(0,H) \subset \mathbb{R}^2$
illustrated in 
Fig.~\ref{fig:symmetry}, in such a way that the air 
temperature $T_{a}$ and panel temperature $T_{w}$ satisfy 
in $\Omega$ the following equations which express conservation of energy:
\begin{subequations}\label{eq:2d_air_panel}
\begin{align}
\nabu \cdot (-k_{a}  \nabu T_a  + \rho_{a} {c_p} 
\widetilde{\vector{v}}_a 
T_a) + \widetilde{\textit{h}}_{aw} (T_a - T_w)  =  0, & 
\label{eq:2d_air}\\[4mm]
\nabu \cdot (-k_{w}  \nabu T_w) + \textit{h}_{aw}^\ast 
(T_w - T_a)  + \textit{h}_{wc}^\ast (T_w - T_c)  = 0, & 
\label{eq:2d_panel}
\end{align}
complemented by the boundary conditions:
\begin{align}
T_a = {T_a}^{in} &  \quad y=0, \label{eq:2d_bc_D}\\[3mm]
-k_{a}  \nabu T_a \cdot \vector{n}  = 0 & \quad y = H,
\label{eq:2d_bc_N_1}\\[3mm]
(-k_{a}  \nabu T_a  + \rho_a c_p \widetilde{\vector{v}}_a T_a) \cdot 
\vector{n} = 0 & \quad x=0, \ x=W, \label{eq:2d_bc_N_2}\\[3mm]
-k_{w}  \nabu T_w \cdot \vector{n} = 0 & \quad y=0,\ y = H,\ x=0,\ x=W.
\label{eq:2d_bc_N_3}
\end{align}
\end{subequations} 

The unknown functions $T_a=T_a(x,y)$ and $T_w=T_w(x,y)$ are the air and wall temperature respectively, $c_p$ is the air specific heat capacity at constant pressure and $k_{a}$ and $k_{w}$ are the thermal conductivities of air and panel material ({\it e.g.}, aluminium), respectively. 
The function $T_c=T_c(x,y)$ represents the temperature of the cooling
two-phase fluid in the channel network and is assumed to be a given
datum in the solution of the equation system~\eqref{eq:2d_air_panel}.
The parameters $\widetilde{\textit{h}}_{aw}$ and 
$\textit{h}_{aw}^\ast$ are the heat transfer coefficient 
$\textit{h}_{aw}$ between air and condenser wall divided 
by suitably defined characteristic lengths $\lambda_{1a}$ and 
$\lambda_{1w}$, respectively. Precisely, 
$\lambda_{1a}$ is related to the variation of 
$k_a$ in the direction between air and condenser wall
while $\lambda_{1w}$ is related to the variation of 
$k_w$ in the thickness of the condenser wall.
The quantity $\textit{h}_{wc}^\ast$ is 
the heat exchange coefficient $\textit{h}_{wc}$ between
the fluid and the panel wall divided by $\lambda_{1w}$. 
The vector field
$\widetilde{\vector{v}}_a$ is the air flow velocity multiplied
by the factor $\lambda_{2a}/\lambda_{1a}$ where $\lambda_{2a}$
is another characteristic length related to the formation of 
the thermal boundary layer at the interface between air and panel.

The quantities 
$W$ and $H$ are the size of the panel in the $x$ and $y$ directions, 
respectively, and $\vector{n}$ is the outward unit vector along the external surface \NEW{of $\Omega$.}
It is  important to notice that the physical properties of air, namely $k_{a}$ and $c_p$, as 
well as the panel material properties, {\it e.g.} $k_{w}$, 
\NEW{and, when simulations are carried out in the natural convection regime, also the
magnitude of the air velocity $\widetilde{\vector{v}}_a$},
depend on the temperatures $T_{a}$ and $T_{w}$, 
hence problem~\eqref{eq:2d_air_panel} is nonlinear.
The detailed derivation of~\eqref{eq:2d_air_panel} 
from the corresponding 3D model
is illustrated for convenience in~\ref{sec:3d_2_2d}.

\subsection{Model for the channel subsystem}
\label{sec:channel_domain}
\label{sec:mathmod1d}
\newcommand{\segment}{\boldsymbol{\sigma}}
\newcommand{\segparameter}{s}
\newcommand{\numel}{M}
\newcommand{\numjunct}{N}
\newcommand{\seglen}{L}
\newcommand{\segstart}{I^{-}}
\newcommand{\segend}{I^{+}}
\newcommand{\vecpos}{\mathbf{x}}
\newcommand{\segdir}{\mathbf{d}}

The channel embedded in each panel, where the two-phase coolant flows, 
is modeled as a pipeline 
network~\cite{bandahertyklar,BrouwerGasserHerty,herty,formaggia2012}, 
{\it i.e.}, a set of \NEW{a number $\numel$ of} 1D straight pipe segments $\segment_{j} \subset \Omega,\ j=1, \ldots, \numel$.
Such segments are joined at a set of $\numjunct$ vertices 
$\vecpos_{i} \in \Omega,\ i=1, \ldots, \numjunct$ and \NEW{each is}
parametrized by a (scalar) local coordinate $\segparameter_{j}$ such 
that $0 \leq \segparameter_{j} \leq \seglen_{j}$, $\seglen_{j}$ being the
length of $\segment_{j}$.
\begin{figure}[h!]
\centering
\includegraphics[width=0.9\columnwidth]{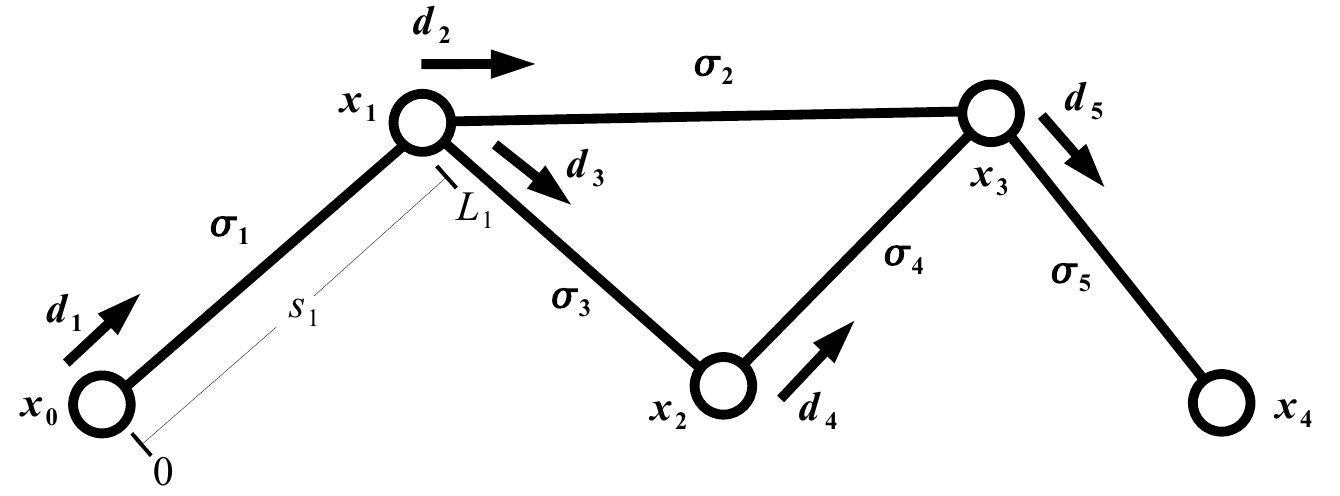}
\caption{Example of channel network geometry and notation.}
\label{fig:channelnotation}
\end{figure}

For each junction $\vecpos_{i}$, we denote by $\segstart_{i} \subseteq 
\{1, \ldots, \numjunct \}$, the set of those indices $j$ for
which $\vecpos_{i}$ is the \emph{first endpoint} 
of the segment $\segment_{j}$, {\it i.e.}, 
$j \in \segstart_{i} \Leftrightarrow \vecpos (\segparameter_{j})\big|_{\segparameter_{j}=0} = \vecpos_{i}$, 
where $\vecpos$ is the (vector) cartesian coordinate.
Similarly, we define $\segend_{i} \subseteq \{1, \ldots, \numjunct \}$, 
to be the set of those indices $j$ for
which $\vecpos_{i}$ is the \emph{second endpoint}
of the segment $\segment_{j}$, {\it i.e.}, 
$j \in \segend_{i} \Leftrightarrow \vecpos (\segparameter_{j})\big|_{\segparameter_{j}=\seglen_{j}} = \vecpos_{i}$.
We assume each parametrization to be \emph{uniform}, {\it i.e.}, 
$j \in \segstart_{m} \cap \segend_{n} 
\Leftrightarrow \vecpos(\segparameter_{j}) = \vecpos_{m} + \segdir_{j} \segparameter_{j}$, where
$\segdir_{j} = (\vecpos_{n}-\vecpos_{m})/\seglen_{j}$ 
is the unit vector defining the \emph{direction}
of $\segment_{j}$.
Furthermore, we introduce the two additional vertices $\vecpos_{0}$ and $\vecpos_{\numjunct+1}$ representing the \emph{inlet}
and the \emph{outlet} of the channel and we assume that they are connected to the first 
node of the first pipe and second node of the last pipe, respectively,
so that we have
$\segstart_{0} \equiv \{ 1 \}$, 
$\segend_{0} \equiv \emptyset$, 
$\segstart_{\numjunct+1} \equiv \emptyset$ and 
$\segend_{\numjunct+1} \equiv \{ \numel \}$.

Within each pipe $\segment_{j}$ we assume the following 1D equations, 
stating conservation of mass, momentum and energy, respectively, to hold 
(see, {\it e.g.}, \cite{hermes2008numerical}):
\begin{subequations}\label{eq:1dns}
\begin{align}
\partial_{\segparameter_{j}} G_{j} = 0, \label{eq:1d_mass}\\
\partial_{\segparameter_{j}} \left( \displaystyle\frac{G_{j}^{2}}{\rho_{j}} + p_{j} \right) =
f_{j} + \rho_{j} \mathbf{g} \cdot \segdir_{j}, \label{eq:1d_momentum}
\\[1mm]
\partial_{\segparameter_{j}} \left( G_{j} \mathcal{H}_{j} \right)+
\textit{h}_{wc} \left( T_{cj} - T_{w} \right) = 0. \label{eq:1d_energy}
\end{align}
\end{subequations}
where $G_{j}$, $\rho_{j}$, $p_{j}$, $f_{j}$, $\mathcal{H}_{j}$ and 
$T_{cj}$
are momentum, density, pressure, frictional forces, specific enthalpy and temperature of the two-phase fluid in each segment $\segment_{j}$, respectively, while $\mathbf{g}$ is the vector denoting the acceleration of gravity.
In view of numerical discretization, it is convenient to rewrite
equations~\eqref{eq:1d_mass}-~\eqref{eq:1d_momentum} as:
\begin{equation}
\label{eq:1dnsplin}
\left\{
\begin{array}{l}
\partial_{\segparameter_{j}} G_{j} = 0\\[2mm]
\partial_{\segparameter_{j}} \varphi_{j}  =
R_{j} G_{j} + \rho_{j} \mathbf{g} \cdot \segdir_{j}
\end{array}
\right.
\end{equation}
where $\varphi_{j} = G_{j}^{2}/\rho_{j} + p_{j}$ denotes the 
total dynamical pressure and $R_{j} = f_{j} / G_{j}$ 
denotes the pipe hydraulic resistance per unit length.
Similarly, equation~\eqref{eq:1d_energy}, upon introducing the 
symbol $\mathcal{W}_{j}$ denoting the enthalpy flux, can be rewritten as:
\begin{equation}
\label{eq:2dnsplin}
\left\{
\begin{array}{l}
\partial_{\segparameter_{j}} \mathcal{W}_{j} = h_{wc} (T_{w} - T_{cj})\\[2mm]
\mathcal{W}_{j} = G_{j} \mathcal{H}_{j}.
\end{array}
\right.
\end{equation}
To close system~\eqref{eq:1dns}, we need:
\begin{enumerate}
\item a set of \emph{coupling} conditions \NEW{at the} $\numjunct$ junctions $\vecpos_{i}, \ i=1, \ldots,\numjunct$;
\item a set of \emph{boundary} conditions \NEW{at the} inlet and outlet
sections; 
\item a set of \emph{constitutive relations}.
\end{enumerate}
All of these relations will be defined in the subsections below.

\subsubsection{Coupling conditions}\label{sec:coupling}
At each of the junction nodes $\vecpos_{i}$
we impose the following \emph{coupling conditions}, 
$\forall i \in \{1, \ldots, \numjunct \},\, \forall j \in \segstart_{i},\, \forall k \in \segend_{i}$:

\begin{subequations}\label{eq:1dnscoupling}

\begin{align}
\left. \varphi_{j} \right|_{\segparameter_{j} = 0} =
\left. \left(\displaystyle\frac{G_{j}^{2}}{\rho_{j}} + p_{j} \right)\right|_{\segparameter_{j} = 0} =
\left. \left(\displaystyle\frac{G_{k}^{2}}{\rho_{k}} + p_{k} \right)\right|_{\segparameter_{k} = \seglen{k}} 
= \left. \varphi_{k} \right|_{\segparameter_{k} = \seglen{k}}, 
\label{eq:cont_tot_press}\\[3mm]
\left. \mathcal{H}_{j}\right|_{\segparameter_{j} = 0} =
\left. \mathcal{H}_{k} \right|_{\segparameter_{k} = \seglen{k}},
\label{eq:cont_enthalpy}\\[4mm]
\displaystyle \sum_{j}\left. -G_{j} \right|_{\segparameter_{j} = 0} +
\displaystyle \sum_{k}\left.  G_{k} \right|_{\segparameter_{k} = 
\seglen{k}} = 0, \label{eq:cons_mass}\\ 
\displaystyle \sum_{j} \left. - \mathcal{W}_{j}
\right|_{\segparameter_{j} = 0} +
\displaystyle \sum_{k} \left. \mathcal{W}_{k}
\right|_{\segparameter_{k} = \seglen{k}} =
\displaystyle \sum_{j} \left.
-\left(G_{j} \mathcal{H}_{j}\right) 
\right|_{\segparameter_{j} = 0} +
\displaystyle \sum_{k} \left.
\left(G_{k} \mathcal{H}_{k}\right) 
\right|_{\segparameter_{k} = \seglen{k}} = 0.
\label{eq:cons_ent_flux}
\end{align}
\end{subequations}
These conditions express continuity of
\emph{total dynamical pressure} and \emph{enthalpy} and 
conservation of mass and enthalpy fluxes at the junctions.

\subsubsection{Boundary conditions}\label{sec:1dnsbc}
At the inlet $\vecpos_{0}$ and outlet $\vecpos_{\numjunct+1}$ we apply the following boundary conditions:
\begin{subequations}\label{eq:1dnscoupling_2}
\begin{align}
\varphi_{1}\Big|_{\segparameter_{1} = 0} = \left(\displaystyle\frac{G_{1}^{2}}{\rho_{1}} + p_{1} \right) \Big|_{\segparameter_{1} = 0} = p_{inlet}, \\
\mathcal{H}_{1} \Big|_{\segparameter_{1} = 0} = 
\mathcal{H}_{inlet}, \\
G_{1} \Big|_{\segparameter_{1} = 0} = 
G_{\numel} \Big|_{\segparameter_{\numel} = \seglen_{\numel}} = G_{tot},
\end{align}
\end{subequations}
where $p_{inlet}$, $\mathcal{H}_{inlet}$ and $G_{tot}$ are 
given data.

\subsubsection{Constitutive relations}\label{sec:1dnsconstrel}
Within each pipe $\segment_{j}$ we assume the following \emph{constitutive relations}, defining the 
\emph{homogeneous flow} regime, to hold:
\begin{subequations}\label{eq:tp_rel}
\begin{align}
\rho_{j} = \displaystyle \frac {\rho_V(T_{cj}) \rho_L(T_{cj})}
{\rho_V(T_{cj}) (1-\textit{x}_{j}) + \rho_L(T_{cj}) \textit{x}_{j}}, \label{eq:tp_rel_rho} \\[4mm]
\mathcal{H}_{j} = \mathcal{H}_L(T_{cj}) 
(1-\textit{x}_{j}) + \mathcal{H}_V(T_{cj}) 
\textit{x}_{j}, \label{eq:tp_rel_h}\\[3mm]
p_{j }= \textit{p}(T_{cj}). \label{eq:tp_rel_p}
\end{align}
\end{subequations}
The two-phase density $\rho_{j}$ and 
enthalpy $\mathcal{H}_{j}$ are calculated through the 
empirical interpolation between all liquid flow (subscript $L$ 
in~\eqref{eq:tp_rel}) 
and all vapor flow quantities (subscript $V$ 
in~\eqref{eq:tp_rel}) that 
are weighted by the vapor quality $\textit{x}_{j}$. 
All the single-phase quantities depend implicitly on the temperature $T_{cj}$
of the two-phase fluid in the $j$-th segment $\segment_{j}$, 
hence system~\eqref{eq:tp_rel} is nonlinear. 
For a detailed description of the two-phase 
constitutive relations, we refer to \cite{collier1996convective}, 
\cite{thome2006engineering} and \cite{whalley1996two}.
After analyzing the review of the most recent correlation 
of the heat transfer coefficient 
for condensation inside tubes \cite{cavallini2003condensation} and \cite{garcia2003review}, 
we have decided to consider 
the \textit{Shah} correlation \cite{shah1979general}, \NEW{valid for 
film condensation pattern, in} the modified version proposed 
in~\cite[Chap. 4]{carichino2010}.
To model the frictional forces $f_{j}$ we used a relation based on the \textit{Blasius equation} \cite[Chap. 13]{thome2006engineering}.
\NEW{For the dependence of the air velocity $\widetilde{\vector{v}}_a$ on the average temperature of the neighbouring panels and on the temperature of air 
at inflow we used correlations given in~\cite{bar1984thermally}.}

\section{Iterative Algorithms}
\label{sec:iter}
\begin{figure}[h!]
\begin{center}
\includegraphics[height=.45\textheight]{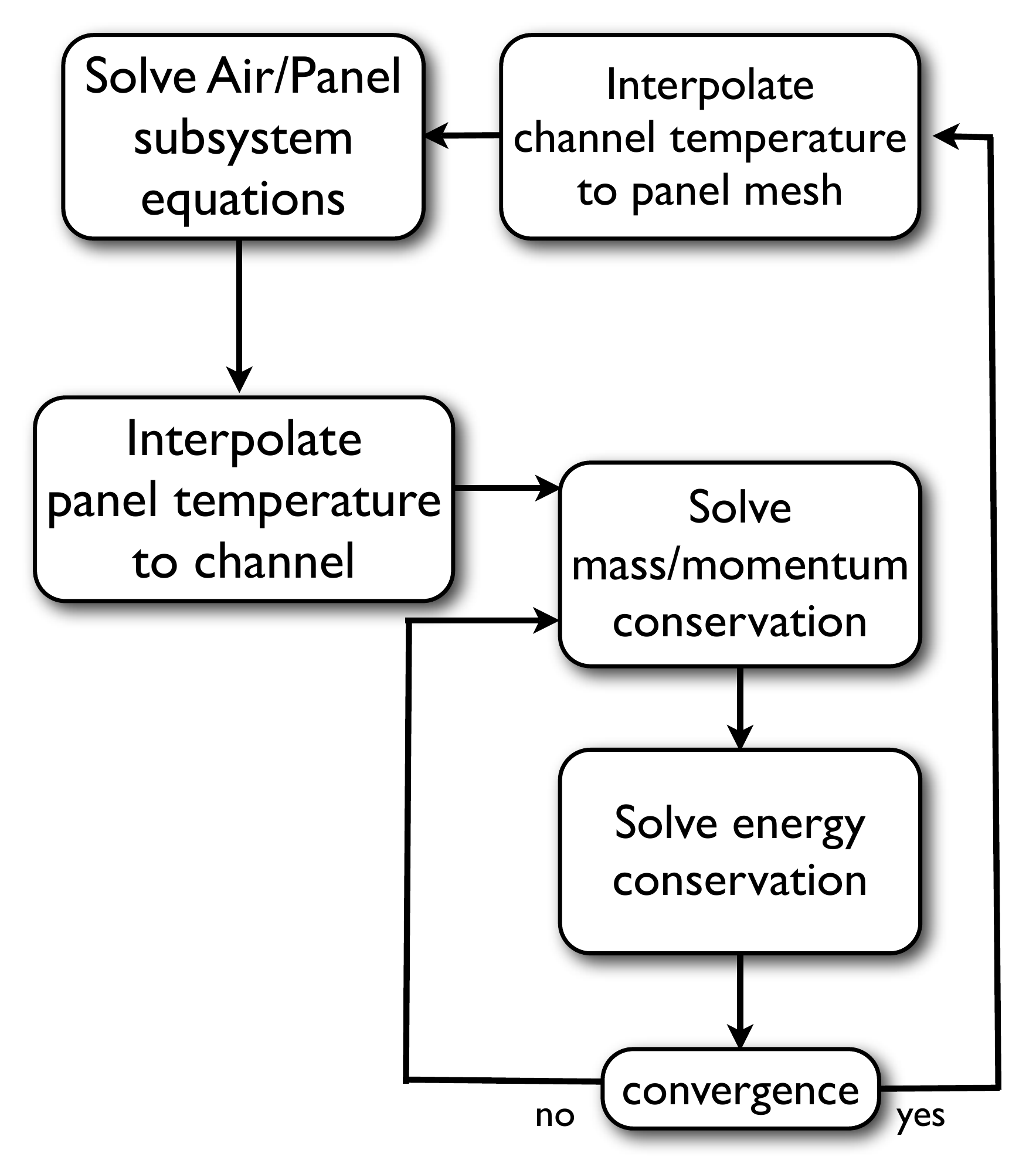}
\caption{Schematic rapresentation of the staggered iteration algorithm.}
\label{fig:iteration}
\end{center}
\end{figure}

The staggered algorithm used for the coupling of the different subsystems is depicted in the flow-chart of Fig.~\ref{fig:iteration}.

The procedure consists of a nested fixed-point iteration \NEW{composed of:
1) an outer iteration loop to solve the 2D air/panel
subsystem; and 2) an inner iteration loop to
solve the non-linear problems within each subsystem.}
In more detail, the outer iteration proceeds as follows:
\begin{enumerate}
\item[1.] Given $T_{c}$ and $\widetilde h_{aw}$, 
compute a new value for $T_{a}$ and $T_{w}$ by solving system~\eqref{eq:2d_air_panel}.
\item[2.] Determine the value of $T_{w}$ at each node of the channel network.
\item[3.] Solve system~\eqref{eq:1dns} to update the quantities $G$, 
$\mathcal{H}$, $p$ describing the state of the fluid
flow in the channel and, as a by-product, 
compute the fluid temperature $T_{c}$ and the vapor quality.
\item[4.] Go back to step 1.
\end{enumerate}

Inner iteration loops are required to solve the non-linear heat flow equations at step 1. and for solving the nonlinear coupled system for the two phase flow at step 3.
For the former we employ a monolithic  quasi-Newton algorithm, while  for the latter we further decouple the equations and proceed as follows:

\begin{enumerate}
\item[3.1] Solve subsystem~\eqref{eq:1d_mass}-~\eqref{eq:1d_momentum}
to update $G$ and $p$.
\item[3.2] Solve subsystem~\eqref{eq:1d_energy}
to update $\mathcal{H}$.
\item[3.3] Determine the density, temperature and vapor quality using 
system~\eqref{eq:tp_rel}.
\item[3.4] Go back to step 3.1.
\end{enumerate}

\section{Dual Mixed-Finite Volume Discretization of the 2D Subproblems} 
\label{sec:air_wall_discr}
In this section we describe the dual mixed-finite volume (MFV)
method used for the numerical approximation 
of the air/panel physical model, presented in 
Sect.~\ref{sec:mathmod2d}.

\subsection{Dual mixed finite element approximation}\label{sec:dm}
Consistently with Sect.~\ref{sec:mathmod2d}, we 
assume that the computational domain $\Omega$ is a 
rectangular open bounded set of $\mathbb{R}^2$ and denote by
$\Gamma:=\partial \Omega$ and $\vector{n}$ 
the domain boundary and its outward unit normal vector, 
respectively. Then, we consider the following 
advection-diffusion-reaction model problem in mixed form: \\
find $u:\Omega \rightarrow \mathbb{R}$ and $\vector{J}: 
\Omega \rightarrow \mathbb{R}^2$ such that: 
\begin{subequations}\label{eq:adr_cf_lp}
\begin{align}
a \vector{J} + \nabu u - a \vector{\beta} u = \vector{0}
& \qquad \mbox{in } \Omega, \label{eq:adr_flux}\\
\nabla \cdot \vector{J} + \gamma u = f & \qquad \mbox{in } \Omega,
\label{eq:adr_balance}\\
u = 0 & \qquad \mbox{on } \Gamma. \label{eq:adr_bcs}
\end{align}
\end{subequations}
In~\eqref{eq:adr_flux}, $a := \alpha^{-1}$ is the inverse
diffusion coefficient, with $\alpha>0$, and the convective 
field $\vector{\beta}$ is a given constant vector, 
while $f \in L^2(\Omega)$ and $\gamma \in L^{\infty}(\Omega)$ 
are given functions, with $\gamma(\vector{x}) > 0$.
We also let $\vector{J}:=-\alpha \nabu u + \vector{\beta} u$ 
be the flux associated with $u$, and we assume 
that (see~\cite{douglasroberts1985,Arbogast_Chen95})
\begin{equation}\label{eq:coercivity}
\dfrac{ \norm{\vector{\beta}}_{L^{\infty}(\Omega)} }
{4 \, \alpha \, \dinf_{\Omega}(\gamma)} < 1.
\end{equation}
Eq.~\eqref{eq:2d_air} 
is a special case of~\eqref{eq:adr_cf_lp} upon setting
$u:=T_a$, $\alpha:=k_a$, $\vector{\beta}:=\rho c_p 
\widetilde{\vector{v}}_a$,
$\gamma:=\widetilde{h}_{aw}$ and $f:= \widetilde{h}_{aw} T_w$, with
$T_w$ a known function, while Eq.~\eqref{eq:2d_panel}
is a special case of~\eqref{eq:adr_cf_lp} upon setting
$u:=T_w$, $\alpha:=k_w$, $\vector{\beta}:=\vector{0}$,
$\gamma:=h_{aw}^\ast + h_{wc}^\ast$ 
and $f:= h_{aw}^\ast T_a + h_{wc}^\ast T_c$, with
$T_a$ and $T_c$ known functions.
Homogeneous Dirichlet boundary conditions \NEW{for $u$} are assumed 
only for ease of presentation, because mixed and/or Neumann
conditions can be easily handled by the proposed scheme 
(see~\cite{sacco1997mixed,BrezziMariniMichPietraSacco2006}).

In view of the numerical approximation of~(\ref{eq:adr_cf_lp}), we introduce 
a regular decomposition $\mathcal{T}_h$ of $\Omega$ into 
${\tt Nel}$ rectangles $K$ of area $|K|$ 
and center of gravity $\vector{x}_{G,K}$,
and we denote by $\mathcal{E}_h$ the set of edges of $\mathcal{T}_h$ 
and by ${\tt Ned}$ the number of total edges of the mesh. We also
let $\mathcal{E}_{h}^{in}$ denote the set of internal
edges of $\mathcal{E}_{h}$. Let $\mathbb{P}_{k_1,k_2}$
be the space of polynomials 
of degree less than or equal to $k_1$ with respect to $x$ and less 
than or equal to $k_2$ with respect to $y$. 
Let $k\geq 0$; for each $K \in \mathcal{T}_h$ we denote 
by $\mathbb{RT}_{[k]}(K) := \mathbb{P}_{k+1,k} (K) \times 
\mathbb{P}_{k,k+1} (K)$ the $k$-th order \textit{Raviart-Thomas} (RT) 
mixed finite element space~\cite{RaviartThomas1977} and 
by $\mathbb{Q}_k(K) = \mathbb{P}_{k,k} (K)$.
We introduce the functional spaces $\mathbf{V} \equiv H_{\Div}(\Omega) = 
\left \{ \vector{v} : \vector{v} \in [L^2(\Omega)]^2, \, 
\nabla \cdot \vector{v} \in L^2(\Omega) \right \}$ and $\mathbf{Q} 
\equiv L^2(\Omega)$, and
their corresponding finite dimensional approximations:
$$
\begin{array}{rcl}
\mathbf{V}_h &=& \left \{ \vector{v}_h \in \mathbf{V} : \  
\vector{v}_h \big|_{K} \in \mathbb{RT}_{[0]}(K) \ \forall K 
\in \mathcal{T}_h \right \},\\
\mathbf{Q}_h &=& \left \{ q_h \in \mathbf{Q} : 
\ q_h \big|_{K} \in \mathbb{Q}_{0}(K) \ \forall K \in 
\mathcal{T}_h \right \}.
\end{array}
$$
Functions in $\mathbf{V}_h$ are linear along each coordinate
direction and discontinuous over $\mathcal{T}_h$ but 
have continuous normal component across each 
edge $e \in \mathcal{E}_{h}^{in}$. Functions
in $Q_h$ are piecewise constant and discontinuous over $\mathcal{T}_h$.

\NEW{To reflect the different nature of
the degrees of freedom of functions in $\mathbf{V}_h$
and $\mathbf{Q}_h$, we introduce two different 
adjacency structures.}

For each (oriented) edge $\vector{e} \in \mathcal{E}_{h}^{in}$, we 
indicate by $|\vector{e}|$ the length of $\vector{e}$, and 
we denote by $K^+_e$ and $K^-_e$ the pair of mesh elements such that
$\vector{e} =\partial K^+_e \cap \partial K^-_e$. We also denote by 
$\vector{n}_e^{+}$ the unit normal vector on $\vector{e}$ 
pointing from $K^+_e$ to $K^-_e$ and define 
$\vector{n}_e^{-} = -\vector{n}_e^{+}$ as the unit normal vector to
$\vector{e}$ pointing from $K^-_e$ to $K^+_e$. In the case where 
$\vector{e} \in \partial \Omega$, we set $\vector{n}_e^{+}: = \vector{n}$.
We indicate by $d_e$ the distance between 
$\vector{x}_{G,K^+_e}$ and $\vector{x}_{G,K^-_e}$.
In the case where $\vector{e} \in \Gamma$, $d_e$ 
is the distance between 
$\vector{x}_{G,K^+_e}$ and the midpoint of edge $\vector{e}$.

For each element $K \in \mathcal{T}_h$, 
we denote by $e(l)$, $l=1, \ldots, 4$, the label number of edge 
$\vector{e}_l$, and by $K_l$ the mesh element neighbour of $K$ with respect
to edge $\vector{e}_l$, whenever $\vector{e}_l$ does not
belong to $\Gamma$.
For any function $w_h \in \mathbf{Q}_h$, 
we introduce the two following operators associated with each edge
of $\mathcal{E}_{h}^{in}$
$$
\llbracket w_h \rrbracket_e := 
w^{K^+_e}\vector{n}_e^+  + w^{K^-_e} \vector{n}_e^{-}, 
\qquad 
\{ w_h \}_e := \Frac{1}{2} (w^{K^+_e} + w^{K^-_e}),
$$
where for each $K \in \mathcal{T}_h$, 
$w^K$ is the constant value of $w_h$ over $K$.
The operator $\llbracket w_h \rrbracket_e$ is the {\em jump} of
$w_h$ across $\vector{e}$ while $\{ w_h \}_l$ is the {\em average}
of $w_h$ across $\vector{e}$. The previous definitions apply also
in the case where $\vector{e} \in \partial{\Omega}$ by setting 
$w^{K^-_e}:=0$. Finally, let $\vector{v}, \, \vector{w}$
be any pair of vectors in $(L^2(\Omega))^2$, and 
\NEW{$v, \, w$} be any function pair in $L^2(\Omega)$. We set 
$A(\vector{v}, \, \vector{w}):= \int_{\Omega} 
a \, \vector{v} \cdot \vector{w}$, $B(v, \, \vector{v}):=
-\int_{\Omega} v \, \nabla \cdot \vector{v}$, $C(v, \, \vector{v}):=
-\int_{\Omega} v \vector{\beta} \cdot \vector{v}$ and
$(v, \, w):=\int_{\Omega} v \, w$.

Then, the dual mixed finite element approximation 
of~(\ref{eq:adr_cf_lp}) over quadrilateral grids reads: 
find $u_h\in \mathbf{Q}_h$ and $\vector{J}_h \in \mathbf{V}_h$ 
such that, for all $\vector{\tau}_h\in\mathbf{V}_h$ and
for all $q_h \in \mathbf{Q}_h$, we have:
\begin{subequations}\label{eq:dm}
\begin{align}
A(\vector{J}_h, \vector{\tau}_h) +
B(u_h, \vector{\tau}_h) +C(u_h, \vector{\tau}_h) =  0 
\label{eq:dm_1}\\[4mm]
B(q_h,  \vector{J}_h) - (q_h, \gamma \, u_h) 
=  -(q_h, f).\label{eq:dm_2}
\end{align} 
\end{subequations}
Equation~\eqref{eq:dm_1} is the discretized 
form of the constitutive law~\eqref{eq:adr_flux}, 
while equation~\eqref{eq:dm_2} is the discretized 
form of the conservation law~\eqref{eq:adr_balance}.
The finite element pair $\mathbf{Q}_h \times \mathbf{V}_h$
satisfies the \textit{inf-sup} compatibility condition, so that
problem~\eqref{eq:dm}, under the coerciveness 
assumption~\eqref{eq:coercivity}, admits a unique 
solution and optimal error estimates can be proved for the 
pair $(u_h, \vector{J}_h)$ in the appropriate graph norm
(see~\cite{RaviartThomas1977,BrezziFortin1991,Arbogast_Chen95}).
The DM formulation can be written in matrix form as
\begin{equation}\label{eq:linear_system_dm}
\left[
\begin{array}{cc}
\mathbf{A} & \quad (\mathbf{B}^T+\mathbf{C}) \\
\mathbf{B} & \quad \mathbf{D} 
\end{array}
\right]
\, 
\left(
\begin{array}{c}
\mathbf{j} \\
\mathbf{u} 
\end{array}
\right)
= 
\left(
\begin{array}{c}
\mathbf{0}_{\tt Ned} \\
\mathbf{f} 
\end{array}
\right)
\end{equation}
where $\mathbf{A} \in \mathbb{R}^{{\tt Ned} \times {\tt Ned}}$ 
is the flux mass matrix,
$\mathbf{B} \in \mathbb{R}^{{\tt Nel} \times {\tt Ned}}$, 
$\mathbf{C} \in \mathbb{R}^{{\tt Ned} \times {\tt Nel}}$ \NEW{and}
$\mathbf{D} \in \mathbb{R}^{{\tt Nel} \times {\tt Nel}}$, 
\NEW{while}
$\mathbf{u} \in \mathbb{R}^{{\tt Nel} \times 1}$, 
$\mathbf{j} \in \mathbb{R}^{{\tt Ned}
\times 1}$ is the unknown vector pair, \NEW{and} $\mathbf{0}_{\tt Ned}$
is the column null vector of size {\tt Ned}.
Two computational difficulties are associated with 
the solution of the DM problem~\eqref{eq:dm}.
The first difficulty is that the linear algebraic 
system~\eqref{eq:linear_system_dm} is in
saddle-point form and has a considerably larger size than
a standard displacement--based method of comparable order.
The second difficulty is that, even in the particular case 
where $\vector{\beta}$ is equal to zero,
it is not possible to ensure that the stiffness matrix acting
on the sole variable $\mathbf{u}$ (obtained upon block 
Gaussian elimination) is an M-matrix for every value of 
$\gamma$ (see~\cite{MariniPietra1989} in the case of triangular RT elements).
This implies that the discrete maximum principle (DMP) can be 
satisfied by the DM method only if the mesh size $h$ is 
sufficiently small, and this constraint may become even more stringent 
if convection is present in the model.

\subsection{The stabilized dual mixed finite volume
approximation}\label{sec:mfv}

To overcome the above mentioned difficulties, we introduce
a (strongly consistent) modification of the DM 
method that extends to the case of quadrilateral grids
the approach \NEW{for triangular grids} proposed  and analyzed 
in~\cite{sacco1997mixed,BrezziMariniMichPietraSacco2006}.
The introduced modifications consist of: 1) replacing the bilinear
form $A(\vector{J}_h, \vector{\tau}_h)$ with the approximate
bilinear form $A_h(\vector{J}_h, \vector{\tau}_h)$ obtained by
using the trapezoidal quadrature formula; 2) replacing 
the bilinear form $C(u_h, \vector{\tau}_h)$ with
$C_h(u_h, \vector{\tau}_h):=C(\{ u_h \}, \vector{\tau}_h)$; 
3) adding to the left-hand side of~\eqref{eq:dm_1} 
the stabilization term
\begin{equation}\label{eq:stabterm}
S(u_h, \, \vector{\tau}_h):= 
- \dsum_{\vector{e} \in \mathcal{E}_h^{in}} \varrho_e(\mathbb{P}e_e)
\dint_{\vector{e}} 
\llbracket u_h \rrbracket_e \cdot \vector{\tau}_h \qquad
\vector{\tau}_h \in \mathbf{V}_h,
\end{equation}
where $\mathbb{P}e_e:= 
(\abs{\vector{\beta} \cdot \vector{n}_e} d_e)/ (2\alpha)$ is the 
{\em local P\`eclet number} associated with edge $\vector{e}$ 
and $\varrho_e: \vector{e} \in \mathcal{E}_h \rightarrow \mathbb{R}^+$ is a 
stabilization function equivalent to adding, for each edge 
of $\mathcal{E}_h$, an {\em artificial diffusion}
to the original problem. 

The resulting stabilized DM formulation reads:
find $u_h\in \mathbf{Q}_h$ and $\vector{J}_h \in \mathbf{V}_h$ 
such that, for all $\vector{\tau}_h\in\mathbf{V}_h$ and
for all $q_h \in \mathbf{Q}_h$, we have:
\begin{subequations}\label{eq:mfv}
\begin{align}
A(\vector{J}_h, \vector{\tau}_h) +
B(u_h, \vector{\tau}_h) + \NEW{C_h(u_h, \vector{\tau}_h)}
+ S(u_h, \, \vector{\tau}_h) =  0 
\label{eq:mfv_1}\\[4mm]
B(q_h, \vector{J}_h) - (q_h, \gamma \, u_h) 
=  -(q_h, f).\label{eq:mfv_2}
\end{align} 
\end{subequations}
The significant advantage of introducing the modifications 1)--3)
with respect to the standard DM approach is that, for each 
element $K \in \mathcal{T}_h$, the flux of $\vector{J}_h$ 
across the edge $e(l)$, $l=1, \ldots, 4$ 
(the degree of freedom of $\vector{J}_h$),
can be expressed {\em explicitly} as a function of the sole
degrees of freedom $u^{K}$ and $u^{K_l}$ as
\begin{equation}\label{eq:approx_flux}
j_{e(l)}(u^K, u^{K}_l) = 
\left[ - \alpha (1 + \varrho_e(\mathbb{P}e_{e(l)}))  
\left( \dfrac{u^{K_l} - u^{K}}{d_{e(l)}} 
\right) + \vector{\beta} \cdot \vector{n}_{e(l)} 
\left( \dfrac{u^{K} +u^{K_l}}{2} \right)\right] |\vector{e}_l|.
\end{equation}
Replacing the above expression into the discrete conservation
law~\eqref{eq:mfv_2}, we end up with the stabilized dual mixed-finite
volume (MFV) approximation of the model problem~\eqref{eq:adr_cf_lp}
\begin{equation}
\label{eq:mfv_stabilized_system}
\left 
\{
\begin{array}{ll}
\dsum_{l=1}^4 j_{e(l)}(u^K, u^{K_l}) + 
u^{K} \overline{\gamma}_{K} |K|
= \overline{f}_{K} |K| & \qquad 
\forall K \in \mathcal{T}_h,\\[5mm]
u^{K_l} = 0 & \qquad \vector{e}_l \in \Gamma,
\end{array} 
\right.
\end{equation}
where $\overline{\gamma}_{K}$ and $\overline{f}_{K}$
are the mean values of $\gamma$ and $f$ on $K$, respectively.
The above proposed stabilized MFV method is the extension to rectangular
elements of the formulation \NEW{for triangular grids} introduced and analyzed 
in~\cite{sacco1997mixed,BrezziMariniMichPietraSacco2006}.
For a similar use of numerical quadrature aimed to construct 
a finite volume variant of the DM method, we refer
to~\cite{vanNooyen} in the case of the 
advection-diffusion-reaction model problem
and to~\cite{Fortin1995} for the 
approximate solution of the Stokes problem in fluid-dynamics.

The MFV method~\eqref{eq:mfv_stabilized_system} 
can be written in matrix form as
\begin{equation}\label{eq:mfv_matrix}
\mathbf{K} \mathbf{u} = \mathbf{g}
\end{equation}
where, for $K=1, \ldots, {\tt Nel}$, 
the entries of the stiffness matrix and of the load vector are:
\begin{equation}\label{eq:mfv_matrix_entries}
\begin{array}{l}
\mathbf{K}_{K,K} = 
\dsum_{l=1}^4
\left[ \dfrac{\alpha (1 + \varrho_e(\mathbb{P}e_{e(l)}))}{d_l}
+ \dfrac{\vector{\beta} \cdot \vector{n}_{e(l)}}{2} 
\right] |\vector{e}_l| + \overline{\gamma}_{K} |K|\\[4mm]
\mathbf{K}_{K,K_l} = 
\left[ -\dfrac{\alpha (1 + \varrho_e(\mathbb{P}e_{e(l)}))}{d_l}
+ \dfrac{\vector{\beta} \cdot \vector{n}_{e(l)}}{2} 
\right] |\vector{e}_l| \\[4mm]
\mathbf{g}_{K} = \overline{f}_{K} |K|.
\end{array}
\end{equation}
Matrix $\mathbf{K}$ is sparse and has at most four nonzero 
entries for each row, in the typical format of lowest-order 
finite volume methods.
Proceeding along the same lines as in~\cite{sacco1997mixed}, 
we can prove the following result.
\begin{proposition}\label{prop:m-matrix}
Let the edge artificial viscosity 
$\varrho_{\vector{e}}(\mathbb{P}e_{\vector{e}})$
be chosen in such a way that
\begin{equation}\label{eq:art_visc}
\varrho_{\vector{e}}(\mathbb{P}e_{\vector{e}}) \geq 
\mathbb{P}e_{\vector{e}} -1 \qquad 
\forall \vector{e} \in \NEW{\mathcal{E}_h.}
\end{equation}
Then, $\mathbf{K}$ is an irreducible diagonally dominant
M-matrix with respect to its colums~\cite{Varga1962}.
\end{proposition}
As a consequence of Prop.~\ref{prop:m-matrix}, the MFV
scheme~\eqref{eq:mfv} satisfies the DMP irrespective 
of the local convective term \emph{and} mesh size.
This lends the scheme a property of robustness
which is a significant benefit in industrial computations like those 
considered in the present article. 
The simplest choice that allows to satisfy~\eqref{eq:art_visc} is
the upwind stabilization
\begin{equation}\label{eq:stab_upwind}
\varrho_l(\mathbb{P}e_{\vector{e}}) = \mathbb{P}e_{\vector{e}}
\qquad \forall \vector{e} \in \mathcal{E}_h.
\end{equation}
Another, more elaborate, choice is the so called Scharfetter-Gummel 
(SG) stabilization
\begin{equation}\label{eq:stab_SG}
\varrho_l(\mathbb{P}e_{\vector{e}}) = 
\mathbb{P}e_{\vector{e}} -1 + \mathcal{B}(2\mathbb{P}e_{\vector{e}})
\qquad \forall \vector{e} \in \mathcal{E}_h,
\end{equation}
where $\mathcal{B}(x):=x/(e^x-1)$ is the inverse of the Bernoulli
function. This latter choice is also known 
as exponential fitting~\cite{ScharfetterGummel1969,RoosStynesTobiska1996}.
The two above stabilizations tend to the same \NEW{limit}
as the P\`eclet number increases. However, their behaviour is quite 
different as the mesh size $h$ decreases, because~\eqref{eq:stab_upwind}
introduces an artificial diffusion of $\mathcal{O}(h)$
as $h \rightarrow 0$ while~\eqref{eq:stab_SG} introduces 
an artificial diffusion of $\mathcal{O}(h^2)$
as $h \rightarrow 0$.
For this reason, the SG stabilized MFV formulation is preferable
as far as accuracy is concerned, and is the one implemented 
in the simulations reported in Sect.~\ref{sec:results}.

\subsection{Numerical validation of the MFV discretization}
\label{sec:numer_validation_MFV}

In this section, we perform a numerical validation of the
stabilized MFV method~\eqref{eq:mfv} applied to the
solution of the model problem~\eqref{eq:adr_cf_lp} with $\Omega=(0,1)\times(0,1)$.

\begin{figure}[h!]
\centering
\subfigure[Upwind stabilization]{
\includegraphics[width=0.75\textwidth]{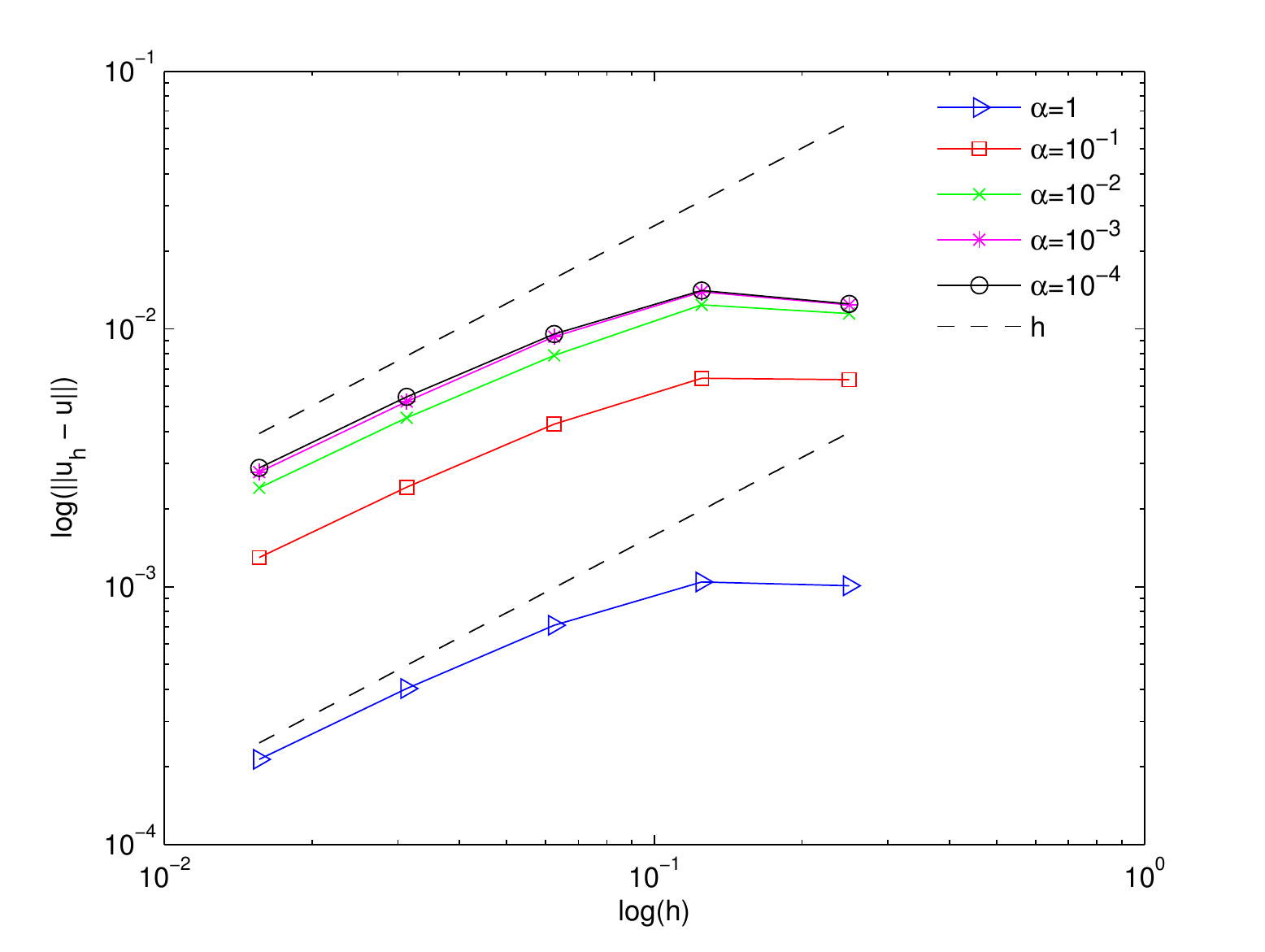}
} 
\subfigure[SG stabilization]{
\includegraphics[width=0.75\textwidth]{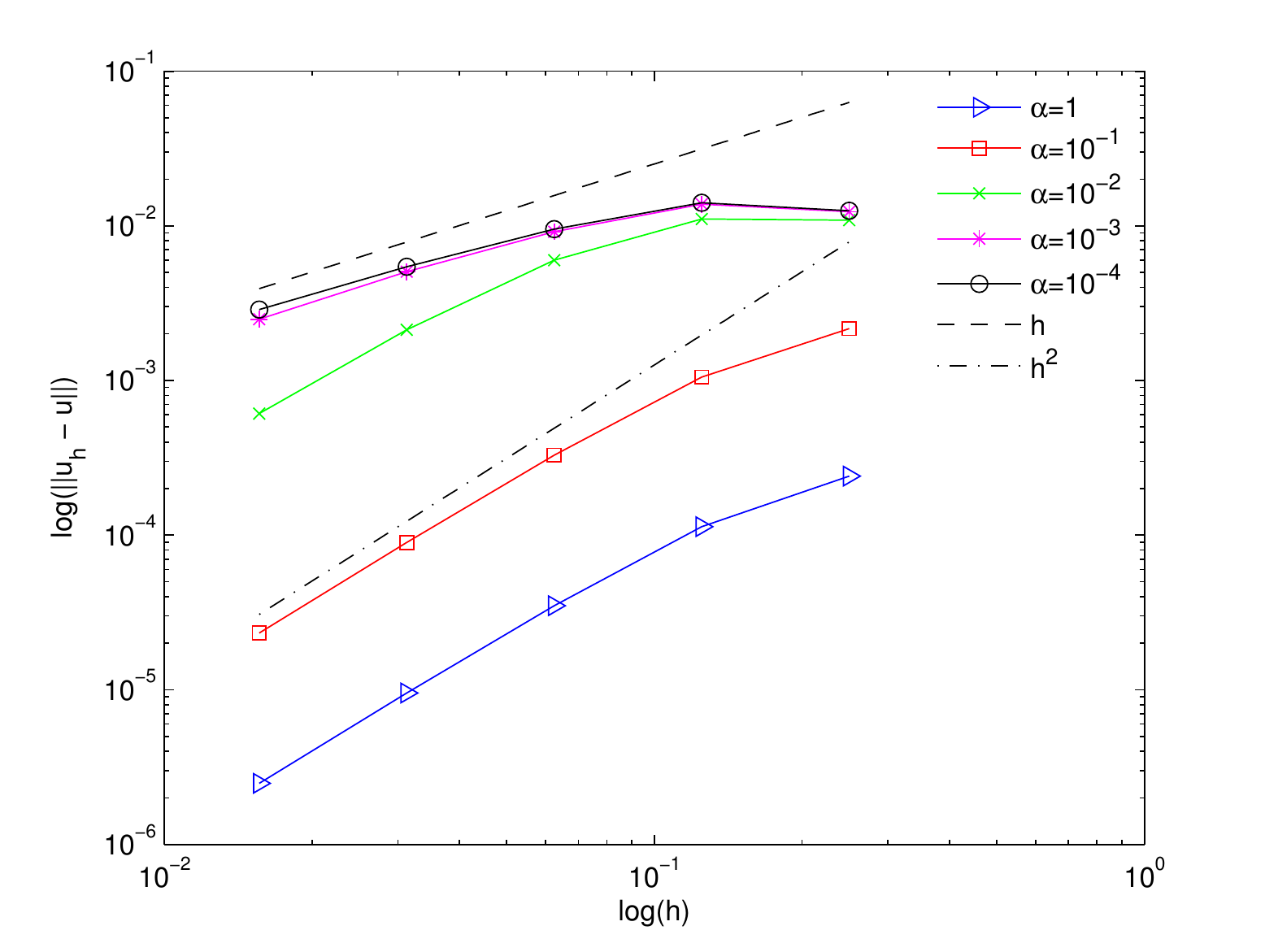}
}
\caption[Estimated convergence rate for upwind and SG methods.]
{Logarithmic plot of the maximum norm of the discretization 
error as a function of $h$ and $\alpha$.}
\label{fig:err}
\end{figure}

In a first case study, we verify the convergence rate
of the scheme when $\vector{\beta}= [0, \ 1]^T$, $\gamma =1$ and 
$f$ is computed in such a way that the exact solution is
$u(x,y) = \cos x \sin y$. As for the diffusion coefficient, we choose
$\alpha= \{ 1,10^{-1},10^{-2},10^{-3},10^{-4} \}$, in order to analyze both dominating diffusive and convective regimes. 
Computations are performed on increasingly refined grids of 
$N \times N$ square elements of dimension varying from $N=4$ to 
$N=64$. Fig.~\ref{fig:err} shows the discrete maximum norm of 
the discretization error 
$$
\| u - u_h \|_{\infty,h}:= \max_{K \in \mathcal{T}_h} 
|u(\vector{x}_{G,K}) - u^K|
$$
as a function of $\alpha$ and of the mesh size $h=1/N$.
Results indicate that for low values of 
the P\`eclet number $\mathbb{P}e$, corresponding for 
example to $\alpha=1$, the SG method has a convergence order of 
$\mathcal{O}(h^2)$, that decreases to $\mathcal{O}(h)$ for 
dominating convection regimes, as for $\alpha=10^{-4}$. 
On the other hand, the estimated convergence error of the 
upwind method is never better than $\mathcal{O}(h)$ for every value
of $\alpha$. 

\begin{figure}[h!]
\centering
\subfigure[Boundary layer example]{
\includegraphics[width=0.75\textwidth]{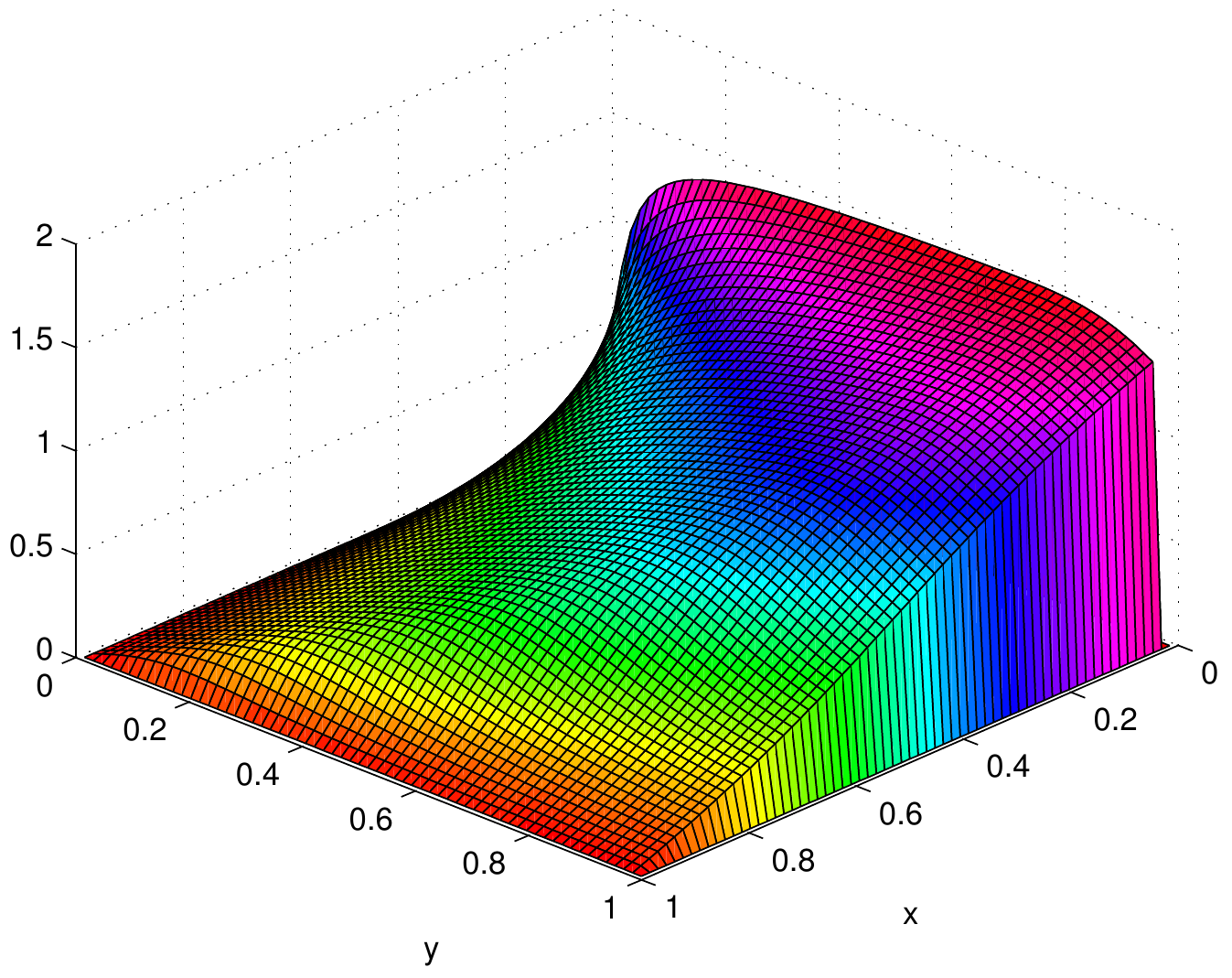}
} 
\subfigure[Interior layer example]{
\includegraphics[width=0.75\textwidth]{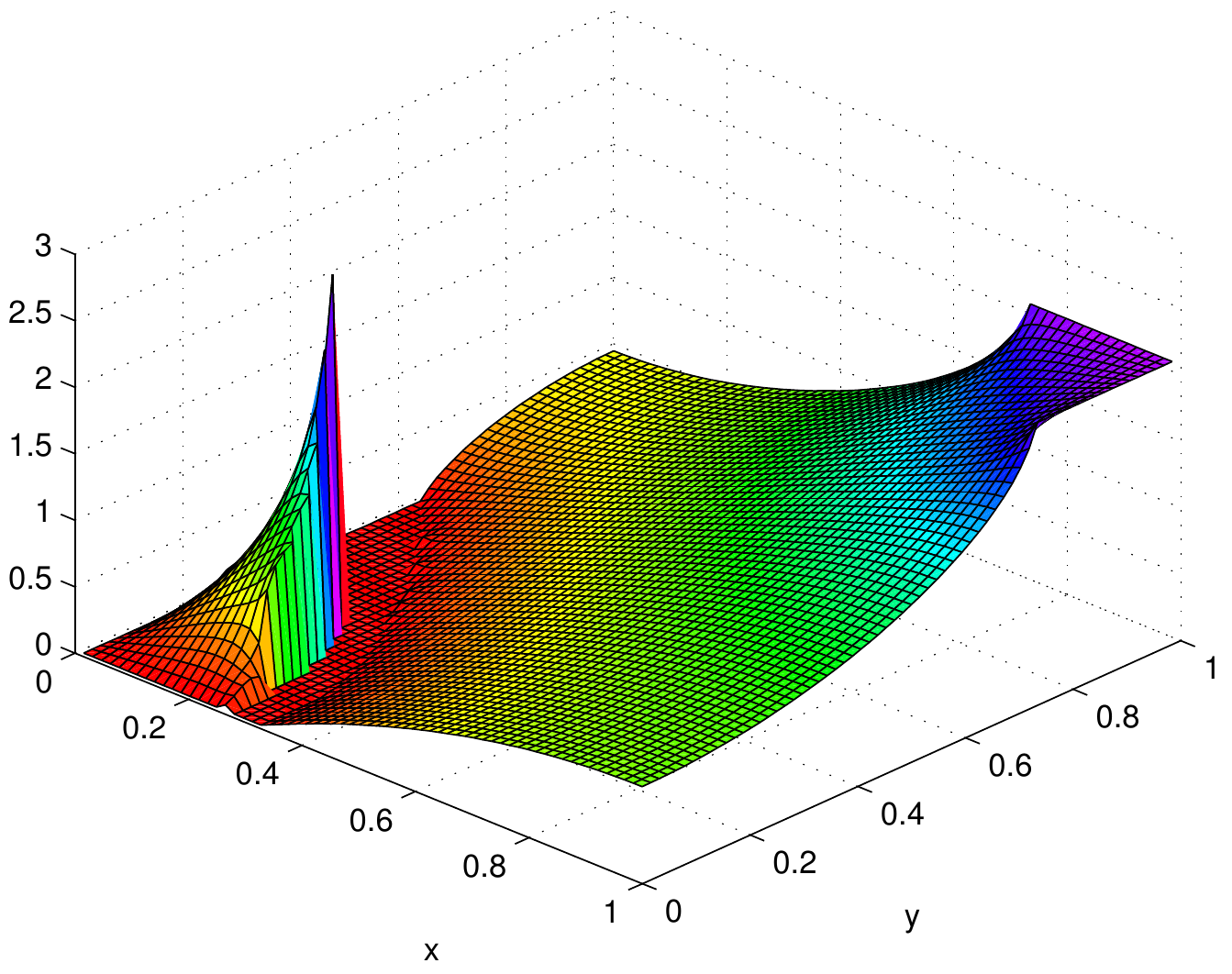}
}
\caption{Surface plot of the numerical solution of two problems with 
sharp boundary and interior layers.}
\label{fig:ex2}
\end{figure}

In a second case study, we validate the robustness and accuracy of the SG 
stabilization in the solution of the two numerical examples 
considered in~\cite{xu1999monotone} where $\alpha = 10^{-6}$, $\gamma=0$
and $h =  2^{-6}$. In the first example, $f=1$
and $\vector{\beta} = [-y, x]^T$. The scope of this computation
is to verify the accuracy and stability of the method in managing 
a boundary layer without introducing spurious oscillations.
In the second example, $f=0$ and $\vector{\beta} = \nabla \psi$,
$\psi$ being the potential function defined as 
$$ 
\psi =
\left \{
\begin{array}{lrcl}
0 & 0 \leq & d + x & < 0.55,\\
2(d - 0.55) & \quad 0.55 \leq & d + x& < 0.65,\\
0.2 & 0.65 \leq &d + x,& 
\end{array} \right.
$$
where $d = ( x^2 + y^2 )^{1/2}$.
Mixed Dirichlet-Robin conditions are enforced on 
the boundary $\Gamma$ in such a way that the solution
exhibits two interior layers, one of which is very sharp.
For graphical purposes, the computed values of $u_h$ have been interpolated through a nodally continuous function.
Results reported in Fig.~\ref{fig:ex2} are in excellent 
agreement with those of~\cite{xu1999monotone} and demonstrate the 
robustness of the stabilized MFV method with respect to 
dominating convective terms and its ability in capturing 
sharp boundary and interior layers without introducing any spurious
oscillation in accordance with
Prop.~\ref{prop:m-matrix}.

\section{Primal Mixed Discretization of the 1D Fluid Equations}
\label{sec:chan_discr}

In this section we focus on the description of the Primal Mixed (PM) finite 
element scheme used for the discretization of the two-phase fluid 
equations~\eqref{eq:1dns} (see~\cite{RobertsThomas1991} for an introduction to
PM methods applied to the numerical solution of elliptic
boundary value problems).
In the following, we consider 
one pipe only and drop the subscript denoting the pipe being considered.

We start by noting that both~\eqref{eq:1dnsplin} and~\eqref{eq:2dnsplin} are special instances
of the following boundary value problem to be solved in the 1D domain $\sigma=(0,L)$:
\begin{subequations}
\label{eq:1dgen}
\begin{align}
\partial_{\segparameter} J = f \label{eq:1dgen_balance} \\
-\varepsilon \partial_{\segparameter} u + \beta u = J + g
\label{eq:1dgen_flux} \\
u(0) = u(L) = 0 \label{eq:1dgen_bcs}
\end{align}
\end{subequations}
where $f$, $\beta$ and $g$ are given data, and $\varepsilon$ is a 
non-negative diffusion coefficient.
We recover~\eqref{eq:1dnsplin} by setting
$J=-G$, $u=\varphi$, $\varepsilon=R^{-1}$, 
$g=-R^{-1} \rho \mathbf{g} \cdot \segdir$, $f=0$ and $\beta=0$,
while we recover~\eqref{eq:2dnsplin} by setting
$J=\mathcal{W}$, $u=\mathcal{H}$, $\varepsilon=0$, 
$g=0$, $f=h_{wc} (T_{w} - T_{c})$ and $\beta=G$.

\begin{remark}[Hyperbolic character of the two-phase fluid model]
From the ma\-the\-ma\-ti\-cal point of view, 
the model problem~\eqref{eq:1dgen}
represents an advective-diffusive model in conservation form
quite similar to that introduced in Sect.~\ref{sec:air_wall_discr} 
for the description of the air/panel physical model.
In the present case, however, there is an important difference because
the two-phase fluid equations~\eqref{eq:1dns} 
have an \emph{hyperbolic} character so that the introduction of 
a diffusive term in the model system~\eqref{eq:1dgen} must be regarded as 
a stabilization term for the corresponding numerical discretization
of equations~\eqref{eq:1dns}. For this reason, throughout the
section, we always assume $\varepsilon$ to be strictly positive.
\end{remark}

\begin{remark}[Extension to pipeline geometry]
The ad\-vec\-tive-dif\-fu\-sive mo\-del \eqref{eq:1dgen} is here solved
in the interval $\sigma=(0,L)$ only for ease of presentation of
the Primal Mixed Finite Element Method approximation. The 
incorporation of the coupling conditions~\eqref{eq:1dnscoupling} 
at each junction node of the pipeline newtwork is straightforward
with the adopted discretization scheme and \NEW{is} discussed
in detail in the remainder of the section.
\end{remark}

Let $\varepsilon$ be a positive bounded function and set 
$a:= \varepsilon^{-1}$. Then, the advective-diffusive 
problem~\eqref{eq:1dgen} can be written in mixed form as: \\
find $u:\sigma \rightarrow \mathbb{R}$ and $J: 
\sigma \rightarrow \mathbb{R}$ such that: 
\begin{subequations}\label{eq:1dgen_mixed}
\begin{align}
a J + \partial_{\segparameter} u - a \beta u + a g = 0
& \qquad \mbox{in } \sigma, \label{eq:1dgen_flux_mixed}\\
\partial_{\segparameter} J = f & \qquad \mbox{in } \sigma,
\label{eq:1dgen_balance_mixed}\\
u(0) = u(L) = 0. & \label{eq:1dgen_bcs_mixed}
\end{align}
\end{subequations}
We assume that
\begin{equation}\label{eq:coercivity_1d}
\partial_{\segparameter}(a \beta) \geq 0 \qquad \mbox{a.e. in } \sigma.
\end{equation}
In view of the numerical approximation of~\eqref{eq:1dgen_mixed}
we introduce a partition $\mathcal{T}_{h}$ of 
$\sigma$ into $N$ intervals $K_{i}$ of 
length $h_i$, $i= 1\ldots N$, 
by means of $N+1$ nodes $s_{j}, j=0\ldots N, s_{0}=0, s_{N}=L$. 
We also introduce the following function spaces defined 
on $\mathcal{T}_{h}$:
\begin{subequations}
\[
\begin{array}{l}
V_{h} = \left\{ v_{h} \in C^{0}(\overline{\sigma}) : 
v_{h}|_{K_{i}} \in \mathbb{P}_{1}(K) \, \forall K \in \mathcal{T}_{h}, \, v_h(0)=v_h(L) = 0\right\}\\[2mm]
Q_{h} = \left\{ p_{h} \in L^{2}(\sigma) : 
p_{h}|_{K_{i}} \in \mathbb{P}_{0}(K) \, \forall K \in \mathcal{T}_{h}\right\}.
\end{array}
\]
\end{subequations}
Functions in $V_h$ are piecewise linear continuous over 
$\overline{\sigma}$ and vanish at the boundary $\partial \sigma$
while functions in $Q_h$ are piecewise constant over $\sigma$.
Nodal continuity of functions in $V_h$ ensures the automatic
satisfaction of the coupling conditions~\eqref{eq:cont_tot_press} 
and~\eqref{eq:cont_enthalpy}.

The PM finite element approximation of~\eqref{eq:1dgen_mixed} 
reads:\\
find $u_h \in V_h$ and $J_h \in Q_h$ such that:
\begin{subequations}\label{eq:PM}
\begin{align}
A(J_h, q_h) + B(u_h, q_h) + C(u_h, q_h) =  -(a g, q_h)  &
\qquad \forall q_h \in Q_h \label{eq:PM_1}\\[4mm]
B(v_h, J_h) =  -(v_h, f) & \qquad \forall v_h \in V_h, 
\label{eq:PM_2}
\end{align} 
\end{subequations}
where 
$$
\begin{array}{l}
A(J_h,q_h):= \displaystyle \int_\sigma a \, J_h \, q_h \, ds, \\[4mm]
B(v_h,J_h):= \displaystyle 
\int_\sigma J_h \, \partial_{\segparameter} v_h \, ds, \\
C(u_h, q_h):= \displaystyle
-\int_\sigma a \, \beta \, u_h \, q_h \, ds
\end{array}
$$
and $(\cdot, \cdot)$ denotes the scalar product in $L^2(\sigma)$.
It can be checked that under the coercivity 
assumption~\eqref{eq:coercivity_1d}, problem~\eqref{eq:PM}
is uniquely solvable.

The PM system~\eqref{eq:PM} can be written in matrix form as
\begin{equation}\label{eq:linear_system_PM}
\left[
\begin{array}{cc}
\mathbf{A} & \quad (\mathbf{B}^T+\mathbf{C}) \\
\mathbf{B} & \mathbf{0}
\end{array}
\right]
\, 
\left(
\begin{array}{c}
\mathbf{j} \\
\mathbf{u} 
\end{array}
\right)
= 
\left(
\begin{array}{c}
\mathbf{g} \\
\mathbf{f} 
\end{array}
\right)
\end{equation}
where $\mathbf{A} \in \mathbb{R}^{N \times N}$ 
is the flux mass matrix,
$\mathbf{B} \in \mathbb{R}^{(N-1) \times N}$ \NEW{and}
$\mathbf{C} \in \mathbb{R}^{N \times (N-1)}$, \NEW{while}
$\mathbf{u} \in \mathbb{R}^{(N-1) \times 1}$,
$\mathbf{j} \in \mathbb{R}^{N \times 1}$ 
is the unknown vector pair, \NEW{and} $\mathbf{0} \in \mathbb{R}^{(N-1) 
\times (N-1)}$ is the null square matrix of size $N-1$.
Compared with the dual mixed system~\eqref{eq:linear_system_dm}, 
the PM formulation~\eqref{eq:linear_system_PM} has a 
considerable advantage because matrix $\mathbf{A}$ is 
\emph{diagonal}, each diagonal entry $A_{kk}$ corresponding to 
the element $K_k$ in the grid, $k=1, \ldots, N$.
Assuming that $\varepsilon$, $\beta$ and $g$ are constant
over each element $K_i$, the first equation 
of~\eqref{eq:PM} can be solved for the flux $J_h$ over each 
mesh element
\begin{equation}\label{eq:flux_PM_K}
J_k = -\varepsilon_k \Frac{u_k - u_{k-1}}{h_k}
+ \beta_k \Frac{u_{k-1}+u_k}{2} - g_k
\qquad \forall i=k, \ldots, N.
\end{equation}
Taking $v_h$ equal to the "hat" function
$\varphi_i$, equal to 1 at every internal node $s_i$ 
and zero at every other node, $i=1, \ldots, N-1$, 
we end up with the following system of \emph{nodal
conservation laws:}
\begin{equation}\label{eq:nodal_system}
J_{i+1} - J_i = f_i \left( \Frac{h_i + h_{i+1}}{2} \right) \qquad
i=1, \ldots, N-1.
\end{equation}
The above equation expresses the fact that at each internal node
of the partition the output flux $J_{i+1}$ is equal to the sum of 
the input flux $J_i$ plus the nodal production term
$P_i:= f_i (h_i + h_{i+1})/2$, in strong analogy with the 
classical Kirchhoff law for the current in an electric circuit.
In particular, if $f=0$, we get strong flux conservation at the
node $s_i$, $i=1, \ldots, N-1$, which corresponds to enforcing
in strong form the coupling conditions~\eqref{eq:cons_mass} and~\eqref{eq:cons_ent_flux}.

Substituting~\eqref{eq:flux_PM_K} into~\eqref{eq:nodal_system} we end
up with the linear algebraic system
in the sole variable $u_h$
\begin{equation}\label{eq:PM_displacement}
\mathbf{M} \mathbf{U} = \mathbf{F}
\end{equation}
where $\mathbf{U} \in \mathbb{R}^{(N-1) \times 1}$ is the
vector of nodal dofs for $u_h$, $\mathbf{F} 
\in \mathbb{R}^{(N-1) \times 1}$ is the right-hand side and 
$\mathbf{M} \in \mathbb{R}^{(N-1) \times (N-1)}$ is the stiffness matrix
whose entries are given by:
$$
M_{ij} = 
\left\{
\begin{array}{ll}
-\Frac{\varepsilon_{i}}{h_{i}} - \Frac{\beta_{i+1}}{2} & 
\qquad j=i-1 \\[2mm]
\Frac{\varepsilon_{i}}{h_{i}} +  
\Frac{\varepsilon_{i+1}}{h_{i+1}} + \Frac{\beta_{i+1}}{2}
-\Frac{\beta_{i}}{2} & \qquad j=i \\[3mm]
-\Frac{\varepsilon_{i+1}}{h_{i+1}} + \Frac{\beta_{i+1}}{2} & 
\qquad j=i+1.
\end{array}
\right.
$$
As in the case of the dual mixed method of Sect.~\ref{sec:dm}, 
the matrix $\mathbf{M}$ turns out to be an M-matrix only
if the mesh size $h$ is sufficiently small.
To avoid this inconvenience, we 
define the local P\`eclet number 
$$
\mathbb{P}e_{i}:= \Frac{|\beta_i| h_i}{2 \varepsilon_i}
\qquad i=1, \ldots, N
$$
and modify the PM finite
element scheme by simply replacing in the first equation of~\eqref{eq:PM}
the term $a=\varepsilon^{-1}$ with 
$$
a_h|_{K_i} := \left(\varepsilon_i(1+ \mathbb{P}e_{i})\right)^{-1}
= \left(\varepsilon_i + \Frac{|\beta_i| h_i}{2} \right)^{-1}
\qquad i=1, \ldots, N.
$$
This amounts to adding a stabilizing artificial diffusion term 
of upwind type (cf.~\eqref{eq:stab_upwind}) into the method and
transforms system~\eqref{eq:PM_displacement} into the stabilized PM scheme
\begin{equation}\label{eq:PM_displacement_stab}
\mathbf{M}^{stab} \mathbf{U} = \mathbf{F}
\end{equation}
where the entries of the stiffness matrix $\mathbf{M}^{stab}$
of the stabilized PM method now read:
$$
M_{ij}^{stab} = 
\left\{
\begin{array}{ll}
-\Frac{\varepsilon_{i}}{h_{i}} - \beta_i^+ & 
\qquad j=i-1 \\[2mm]
\Frac{\varepsilon_{i}}{h_{i}} +  
\Frac{\varepsilon_{i+1}}{h_{i+1}} + \beta_{i+1}^+ - \beta_{i}^- 
& \qquad j=i \\[3mm]
-\Frac{\varepsilon_{i+1}}{h_{i+1}} + \beta_{i+1}^- & 
\qquad j=i+1
\end{array}
\right.
$$
having set:
$$
\begin{array}{ll}
\beta^+:= \Frac{\beta + |\beta|}{2} & \qquad (\geq 0) \\
\beta^-:= \Frac{\beta - |\beta|}{2} & \qquad (\leq 0).
\end{array}
$$
By inspection on the expressions of $M_{ij}^{stab}$ we have the 
following result.
\begin{proposition}\label{prop:m-matrix_PM}
The stiffness matrix $\mathbf{M}^{stab}$ is an irreducible diagonally 
dominant M-matrix with respect to its colums. 
\end{proposition}
As in the case of the MFV scheme, 
Prop.~\ref{prop:m-matrix_PM} implies that the upwind
stabilized PM finite element scheme satisfies the DMP.
Moreover, the upwind PM method is at most 
first-order accurate with respect to the discretization parameter $h$. 
\begin{remark}[Stabilization method]
In the case of problem~\eqref{eq:2dnsplin} 
the SG stabilization~\eqref{eq:stab_SG} cannot be used because 
$\varepsilon=0$. Therefore, to ensure a consistent treatment that is 
applicable in both hyperbolic and advective-diffusive regimes, the 
artificial diffusion term of upwind type~\eqref{eq:stab_upwind} 
is added in the numerical examples of Sect.~\ref{sec:numer_validation_PM}
\NEW{and of Sect.~\ref{sec:results}.}
\end{remark}

\subsection{Numerical validation of the PM discretization}
\label{sec:numer_validation_PM}

\begin{figure}[h!]
	\centering
		\includegraphics[width=\textwidth]{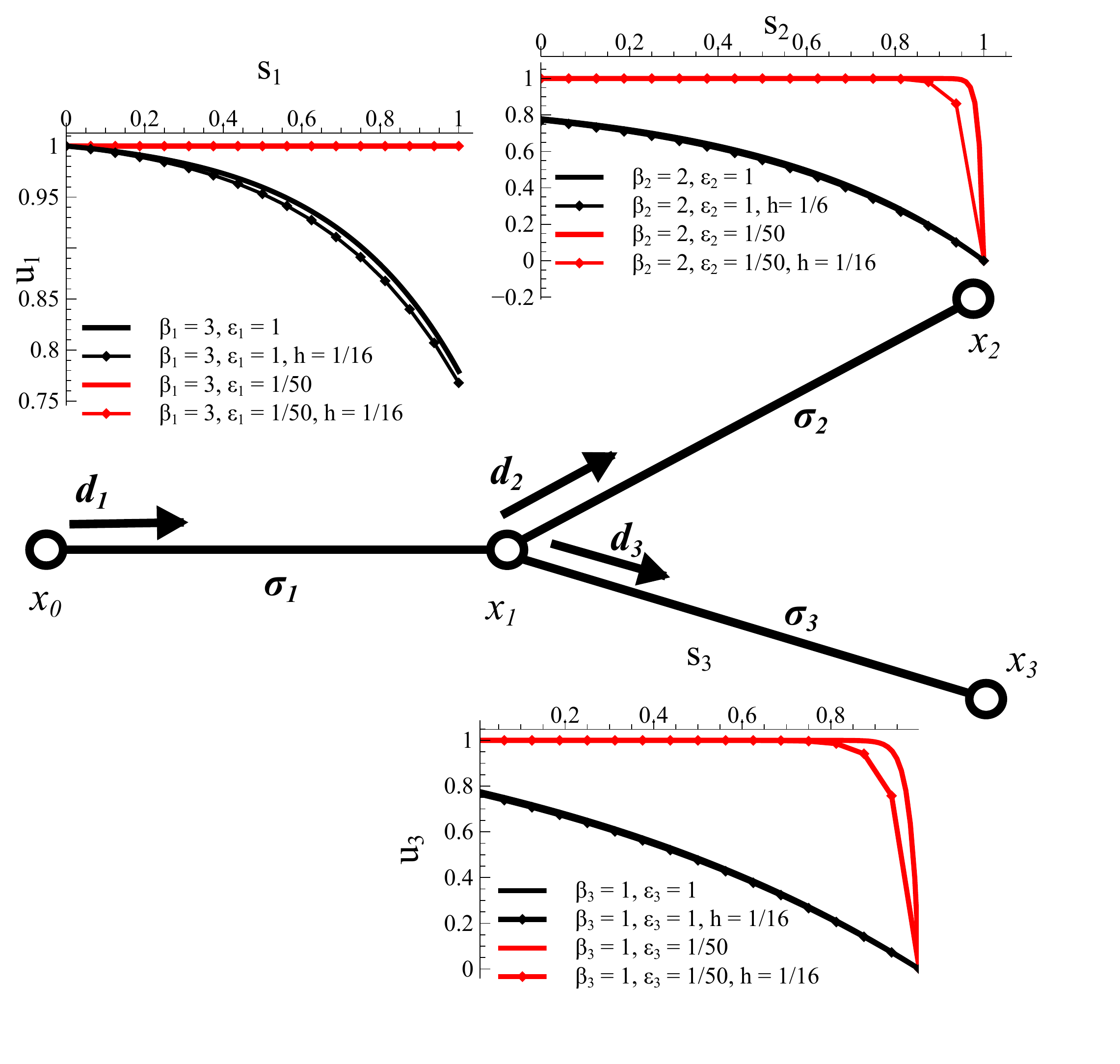}
	\caption{Test case on a three-segment network. 
	Solid lines denote the exact solution while dotted lines denote the 
	numerical solution computed by the upwind stabilized PM method.}
	\label{fig:test_case_network1}
\end{figure}

In this section, we perform a numerical validation of the stabilized PM
method~\eqref{eq:PM} applied to the solution of the model problem~\eqref{eq:1dgen} 
on the test network geometry depicted in Fig.~\ref{fig:test_case_network1}.

In the first test case we study a diffusion-dominated flow while 
in the second test case the flow is in the advection-dominated regime.
For both cases we let $\beta|_{\boldsymbol{\sigma}_{1}} = 3$,
$\beta|_{\boldsymbol{\sigma}_{2}} = 2$, $\beta|_{\boldsymbol{\sigma}_{3}} = 1$ 
and $f = g = 0$.
For the first test case, whose exact solution is shown in black in 
Fig.~\ref{fig:test_case_network1}, we let $\varepsilon = 1$ on all network segments,
while for the second test, whose exact solution is shown in red in 
Fig.~\ref{fig:test_case_network1}, we let $\varepsilon = 1/50$.
It is easily verified that the exact solution of both tests can be expressed
\NEW{as:}
\begin{subequations}\label{eq:exact_solution_PM}
\begin{align}
u(s) |_{\boldsymbol{\sigma}_{i}}&= u(0)|_{\boldsymbol{\sigma}_{i}} \Frac{ e^{\displaystyle (\beta|_{\boldsymbol{\sigma}_{i}} L_{i})/ \varepsilon} -  e^{\displaystyle (\beta|_{\boldsymbol{\sigma}_{i}} s)/ \varepsilon} }{  e^{\displaystyle (\beta|_{\boldsymbol{\sigma}_{i}} L_{i})/ \varepsilon} -1 } + u(L_{i})|_{\boldsymbol{\sigma}_{i}} \Frac{ e^{\displaystyle (\beta|_{\boldsymbol{\sigma}_{i}} s)/ \varepsilon} -  1}{  e^{\displaystyle (\beta|_{\boldsymbol{\sigma}_{i}} L_{i})/ \varepsilon} -1 } 
\label{eq:exact_u_s}\\
J(s) |_{\boldsymbol{\sigma}_{i}}&= \beta|_{\boldsymbol{\sigma}_{i}} \Frac{ u(0)|_{\boldsymbol{\sigma}_{i}} e^{\displaystyle (\beta|_{\boldsymbol{\sigma}_{i}} L_{i})/ \varepsilon} - u(L_{i})|_{\boldsymbol{\sigma}_{i}} }{  e^{\displaystyle (\beta|_{\boldsymbol{\sigma}_{i}} L_{i})/ \varepsilon} -1 }
\label{eq:exact_J_s}
\end{align}
\end{subequations}
for $i=1,2,3$, where $L_1=L_2=L_3=1$, $u(0) |_{\boldsymbol{\sigma}_{1}}=1$, $u(L_1) |_{\boldsymbol{\sigma}_{1}}= u(0) |_{\boldsymbol{\sigma}_{2}} = u(0) |_{\boldsymbol{\sigma}_{3}} = \omega$ and $u(L_2) |_{\boldsymbol{\sigma}_{2}} = u(L_3) |_{\boldsymbol{\sigma}_{3}}=0$. 
The value $\omega$ of the solution $u$ at the junction node $x_1$
is determined from the flux continuity condition
$$
J|_{\boldsymbol{\sigma}_{1}} = J|_{\boldsymbol{\sigma}_{2}} + J|_{\boldsymbol{\sigma}_{3}},
$$
that yields
$$
\omega = \Frac{ \frac{\varepsilon}{L_1} \mathcal{B}\left( \frac{-\beta|_{\boldsymbol{\sigma}_{1}} L_1}{\epsilon} \right) u(0)|_{\boldsymbol{\sigma}_{1}} +
\frac{\varepsilon}{L_2} \mathcal{B}\left( \frac{\beta|_{\boldsymbol{\sigma}_{2}} L_2}{\epsilon} \right) u(L_2)|_{\boldsymbol{\sigma}_{2}} +
\frac{\varepsilon}{L_3} \mathcal{B}\left( \frac{\beta|_{\boldsymbol{\sigma}_{3}} L_3}{\epsilon} \right) u(L_3)|_{\boldsymbol{\sigma}_{3}}
}{
\frac{\varepsilon}{L_1} \mathcal{B}\left( \frac{\beta|_{\boldsymbol{\sigma}_{1}} L_1}{\epsilon} \right) +
\frac{\varepsilon}{L_2} \mathcal{B}\left( \frac{-\beta|_{\boldsymbol{\sigma}_{2}} L_2}{\epsilon} \right)  +
\frac{\varepsilon}{L_3} \mathcal{B}\left( \frac{-\beta|_{\boldsymbol{\sigma}_{3}} L_3}{\epsilon} \right)
 }
$$
where $\mathcal{B}$ is the inverse of the Bernoulli function 
introduced in Sect.~\ref{sec:mfv}.

\begin{figure}[h!]
\centering
\subfigure[$\| u - u_h \|_V$]{
\includegraphics[width=0.45\textwidth]{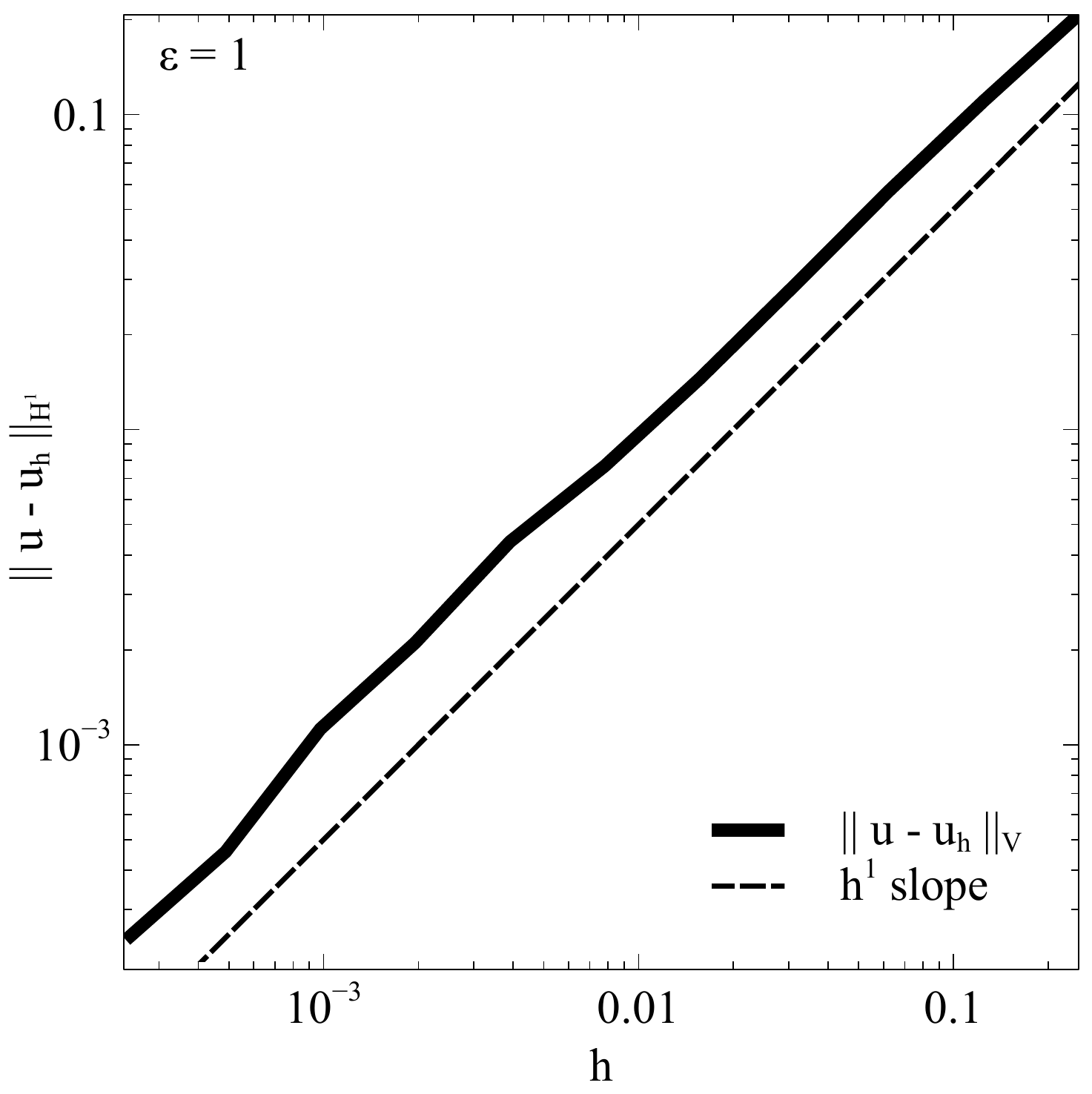}
} 
\subfigure[$\| J - J_h \|_Q$]{
\includegraphics[width=0.45\textwidth]{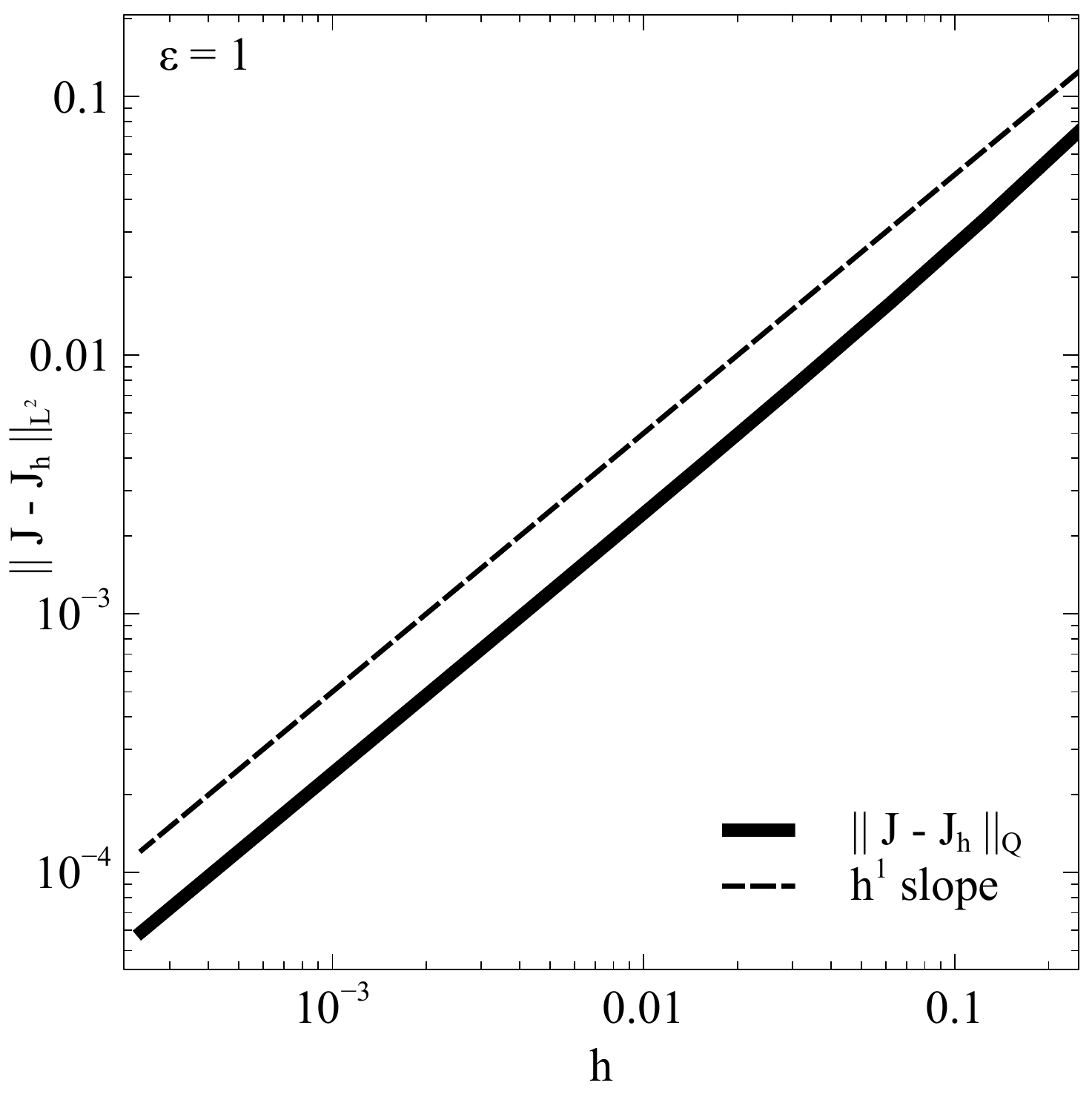}
}
\caption{Logarithmic plot of the discretization error as a 
function of $h$ in the case $\varepsilon=1$.}
\label{fig:err_diffusive_regime}
\end{figure}

Fig.~\ref{fig:err_diffusive_regime} 
shows the logarithmic plots of the discretization errors
$\| u - u_h\|_V$ and $\| J-J_h\|_Q$ as a function of the 
discretization parameter $h$ in the diffusive-dominated regime.
The scheme turns out to have a first-order accuracy. 
This \NEW{result confirms the validity of the error analysis 
carried out in~\cite{RobertsThomas1991} 
in the case of a purely diffusive problem also in the case of an 
advective-diffusive model.}

\begin{figure}[h!]
\centering
\subfigure[$\| u - u_h \|_V$]{
\includegraphics[width=0.45\textwidth]{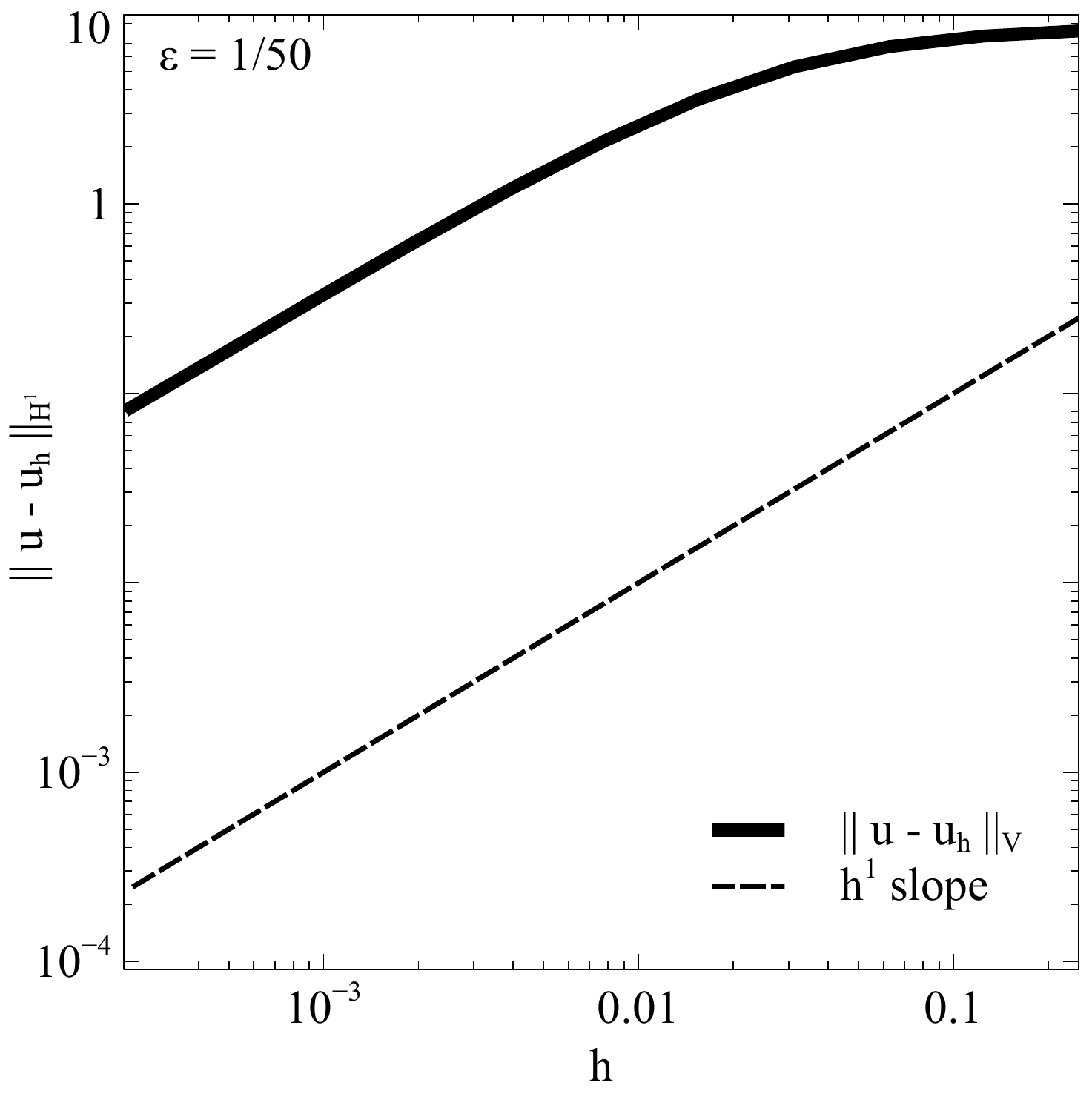}
} 
\subfigure[$\| J - J_h \|_Q$]{
\includegraphics[width=0.45\textwidth]{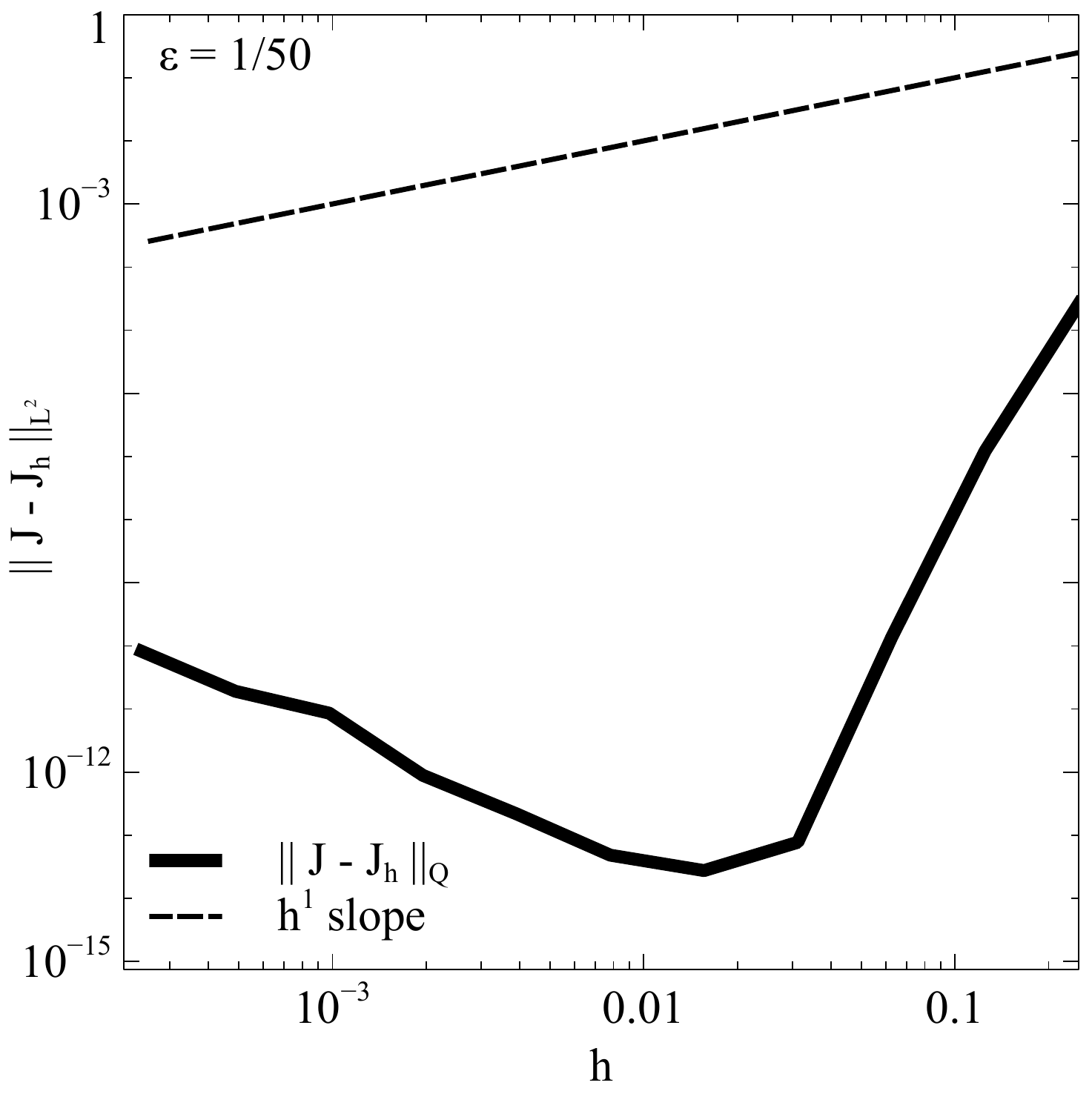}
}
\caption{Logarithmic plot of the discretization error as a 
function of $h$ in the case $\varepsilon=1/50$.}
\label{fig:err_advective_regime}
\end{figure}

Fig.~\ref{fig:err_advective_regime} 
shows the logarithmic plots of the discretization errors
$\| u - u_h\|_V$ and $\| J-J_h\|_Q$ as a function of the 
discretization parameter $h$ in the advective-dominated regime.
The scheme is still first-order accurate 
with respect to $h$ in the computation of the primal
variable $u$ despite the fact that the magnitude of the computed error 
is higher than in the diffusion-dominated regime.
\NEW{The reported error curve for the flux variable $J$ is dominated by the effect of round-off}, 
in accordance with the fact that
in the advective-dominated regime the flow is almost hyperbolic
and the computed flux $J_h$ is a very good approximation of the 
exact flux $J$.

We conclude the validation analysis of the upwind 
stabilized PM method by considering again 
Fig.~\ref{fig:test_case_network1} which shows the numerical 
solution of the benchmark problem (denoted by black and red 
dotted curves) computed with a grid 
spacing $h = 1/16$ and superposed to the exact 
solution~\eqref{eq:exact_u_s}. It is to be noted 
that in the advective-dominated regime ($\varepsilon=1/50$)
the numerical solution almost coincides with the exact one
in the first branch of the network $\sigma_1$ because 
there the problem is almost hyperbolic and the input datum
$u(0)=1$ is transported by the fluid velocity.
We also note that in the other two branches of the network,
$\sigma_2$ and $\sigma_3$, even though the chosen stepsize is not 
sufficiently small to fully resolve the boundary layer
at the outlets, the PM upwind method provides a solution which
is monotone and free of spurious oscillations in accordance with
Prop.~\ref{prop:m-matrix_PM}. A more considerable error occurs
in the computed solution when the problem is diffusion-dominated 
($\varepsilon=1$) in accordance with the fact that the PM is only
first-order accurate.

\section{Simulation Results}\label{sec:results}

In this section we \NEW{perform} a thorough validation of the computational model illustrated in the previous sections. The simulations are representative of realistic geometries of advanced cooling systems for power electronics. In particular aluminum condenser panels, as part of a two-phase thermosyphon loop, are simulated in natural convection operation mode. 
In Sect.~\ref{sec:numres1} we analyze the impact of \NEW{channel geometry and topology} on the cooling performance, while
in Sect.~\ref{sec:numres2} we compare the model predictions with the measured data reported in~\cite{iecon11} and based on the experimental campaign and methodology illustrated in~\cite{experimental}.

\subsection{Comparison of different channel geometries}\label{sec:numres1}
\NEW{
In this section we use our simulation code to estimate the impact of different
pipe geometries on the cooling properties of the system. 
With this aim, we consider three test cases where panel size and material, input power, air velocity 
and temperature are the same, but with different channel paths.}

\NEW{The developed code represents a strong tool in the design of complex channel geometries allowing the researchers to optimize the topology of complex systems.}

\NEW{The simulation data are summarized in Tab.~\ref{tab:tabdata_2}.}

\begin{table}[!h]
\begin{center}
\begin{tabular}{|c|c|c|}
\hline
Parameter & Value & Units \\
\hline $S$ & $0.05$ & m \\
\hline $\lambda_1=\lambda_2$ & $0.025$ & m \\
\hline $T_a^{in}$ & $298.15$ & K \\
\hline $|\vector{V}_a^{in}|$ & $1$ & ${\rm m} \, {\rm s}^{-1}$ \\
\hline $W$ & $0.45$ & m \\
\hline $H$ & $0.2$ & m\\
\hline $T_0$ & $358.15$ & K\\
\hline $G_{tot}$ & $5.8$ & 
${\rm Kg} \, {\rm m}^{-2} \, {\rm s}^{-1}$ \\
\hline $h_{wc}$ & $3$ & ${\rm W} \, {\rm m}^{-2} \, {\rm K}^{-1}$ \\
\hline $h_{aw}$ & $1.1$ & ${\rm W} \, {\rm m}^{-2} \, {\rm K}^{-1}$ \\
\hline
\end{tabular}
\end{center}
\caption{Model parameters.}
\label{tab:tabdata_2}
\end{table}

\begin{figure}[h!]
\centering
\subfigure[device A]{\includegraphics[width=.7\linewidth]{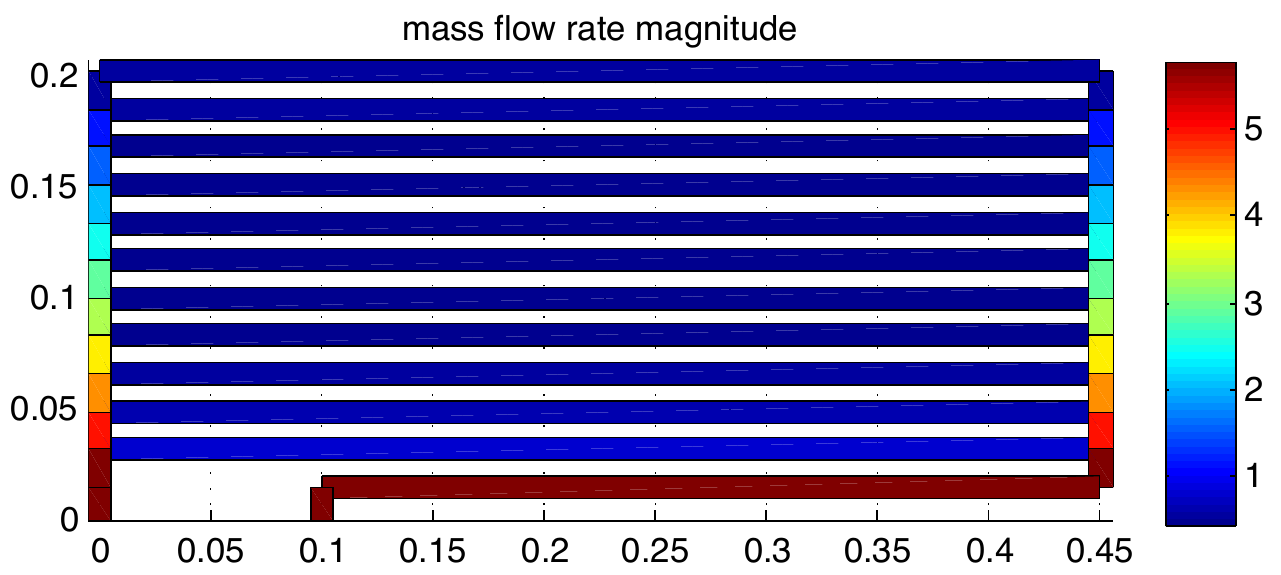}}
\subfigure[device B]{\includegraphics[width=.7\linewidth]{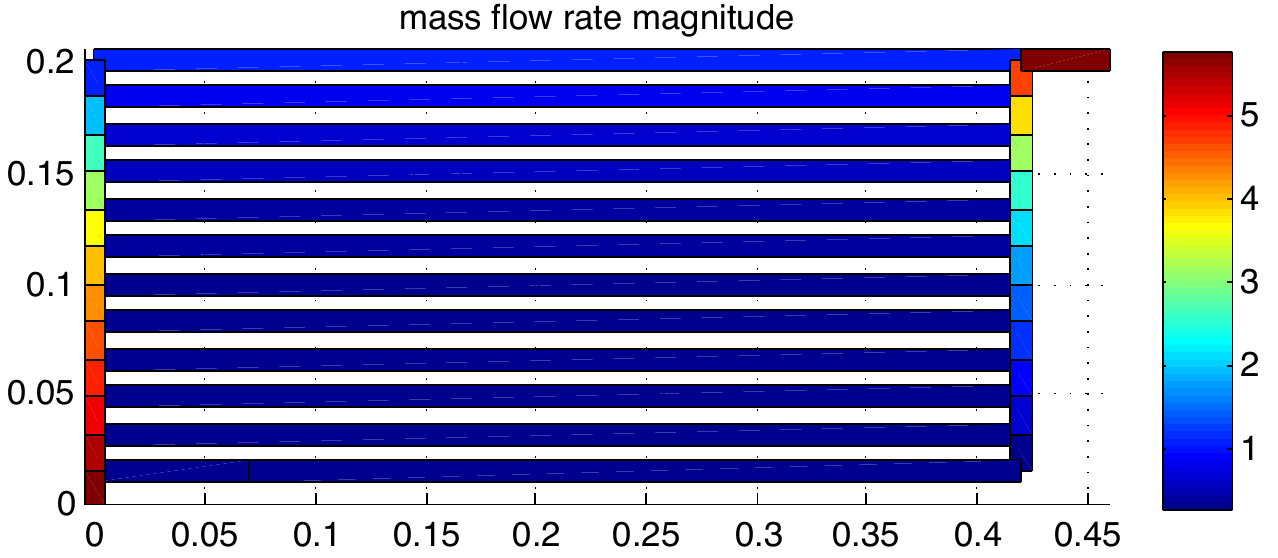}}
\subfigure[device C]{\includegraphics[width=.7\linewidth]{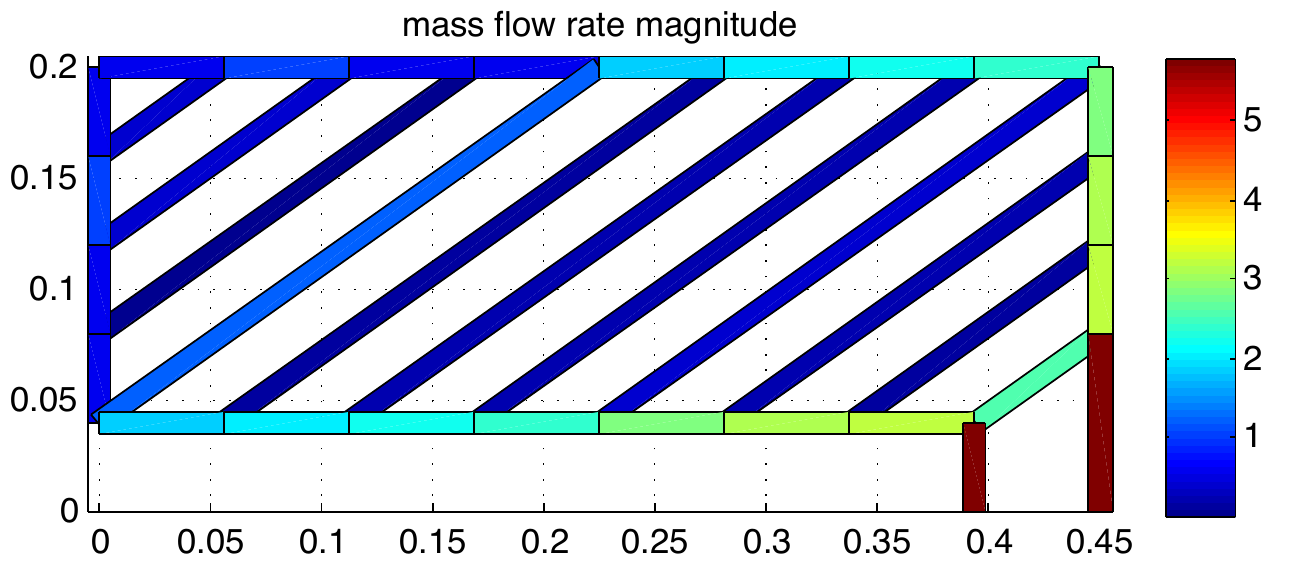}}
\caption{Comparison of mass flow rate magnitude for the three devices geometries.}
\label{fig:geocomp}
\end{figure}

\NEW{The mass flux of coolant is the input datum of the simulation. 
Such value represents the total mass flowing through the panel, assuming 
the coolant to be in full vapor state at the inlet of the system.}

\NEW{The geometry of the three devices is compared in Fig.~\ref{fig:geocomp}, 
with the color scale representing the absolute value of the mass flow rate
in each channel segment.}

\NEW{The structure of a condenser panel is based on a series of parallel channels.
A good flow distribution is a mandatory element for an optimal design, allowing the designer to maximally exploit the
system and therefore increasing the maximum power density of the cooling device. Case "a" and case "b" indicate
a better distribution of mass flow over the parallelized channels compared to case "c". Starting from case "a" and
"b", we see that the flow distribution is a function of the flow-path resistance: the higher the flow-path resistances, the
lower is the flow rate. For case "a", the flow rate is higher in the lower channels, closer to the inlet, and slightly decreases
toward the top part of the panel. The configuration "b" 
is a possible design solution to overcome the pressure drop unbalance that may occur among
the channels, and to guarantee a more uniform distribution over the entire surface due to equal inlet-outlet channel-flow-path length. 
Unfortunately, this effect is not present and a distribution of the mass flow rate similar to that 
in case "a" is obtained. Case "c" is studied to take advantage of the channel orientation 
and the positive effect of the gravitational field in the 
condensation process.}
\begin{figure}[h!]
\centering
\subfigure[device A]{\includegraphics[width=.7\linewidth]{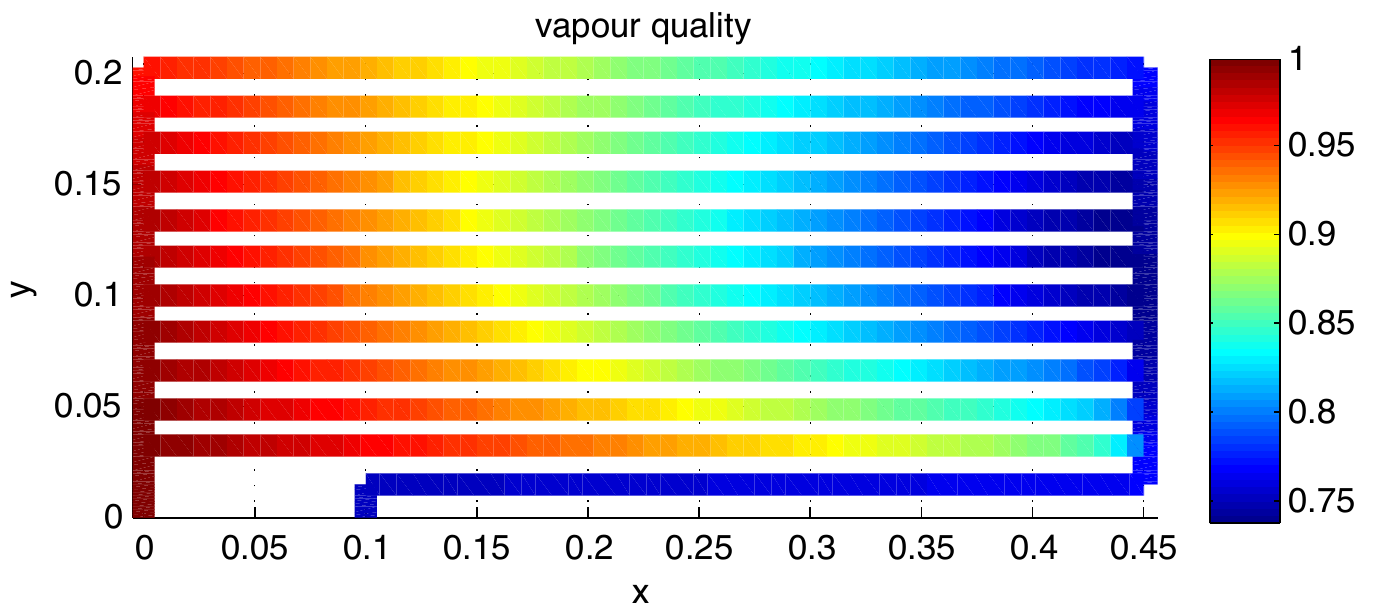}}
\subfigure[device B]{\includegraphics[width=.7\linewidth]{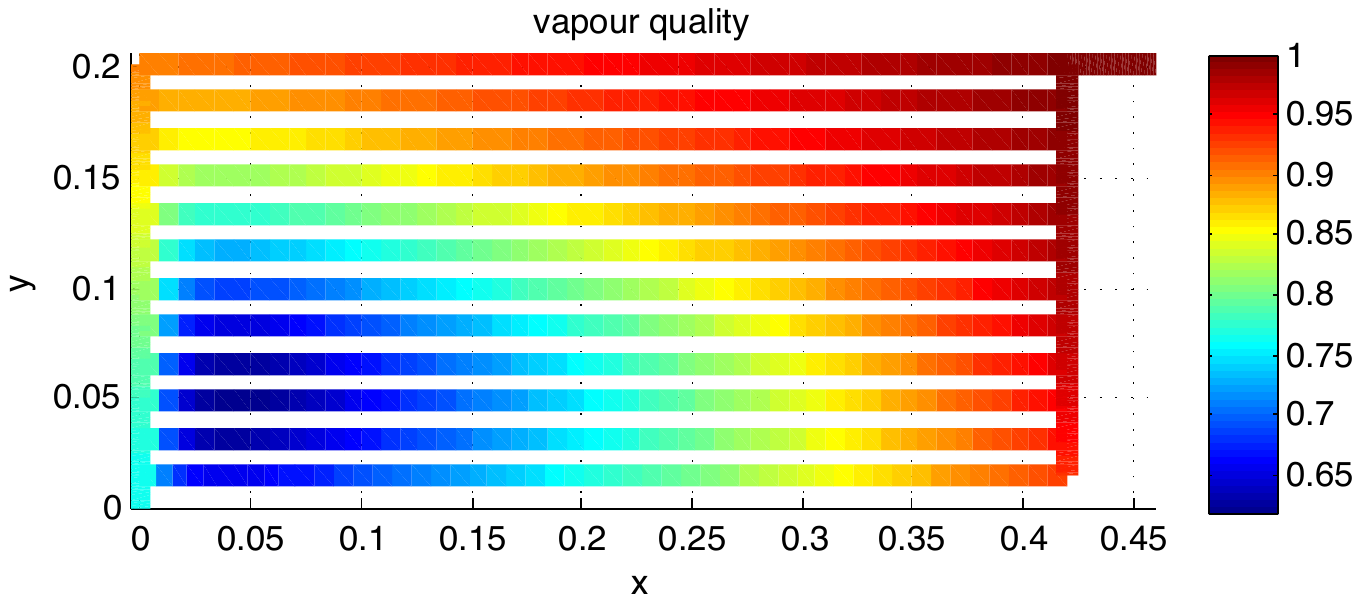}}
\subfigure[device C]{\includegraphics[width=.7\linewidth]{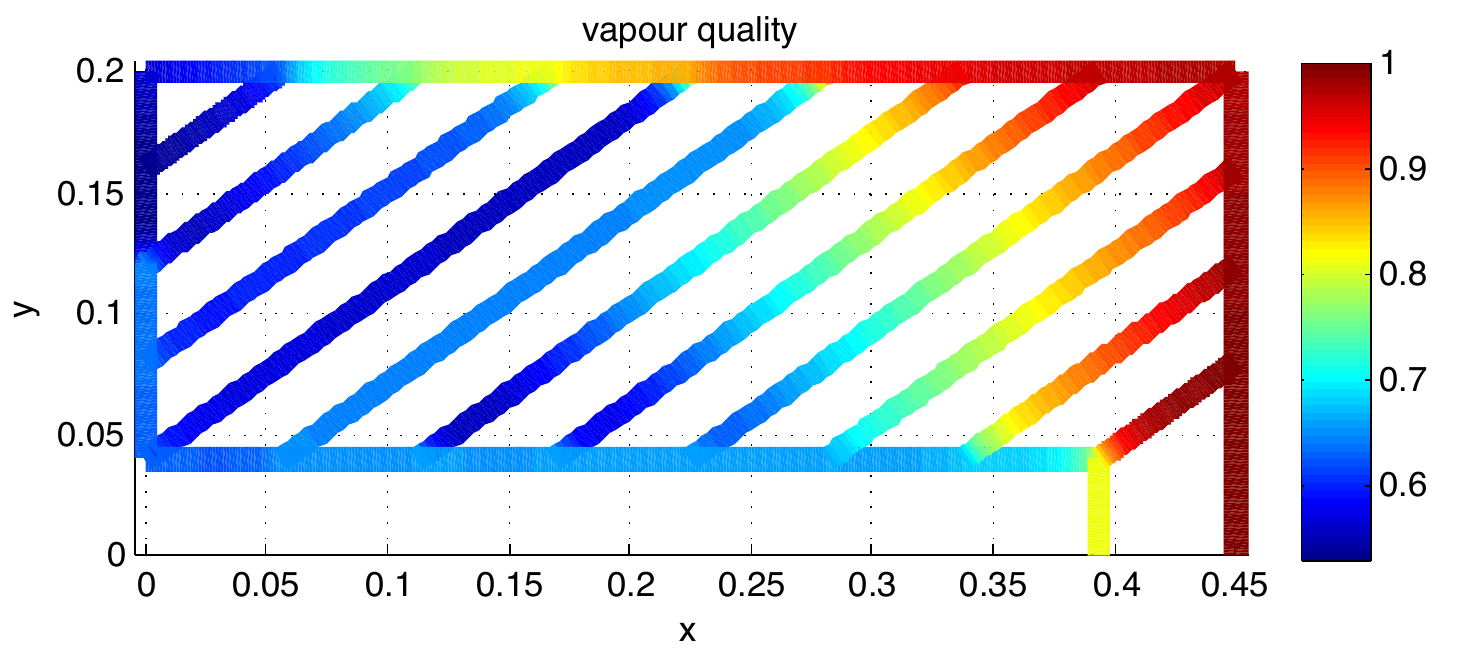}}
\caption{Comparison of vapor quality for the three devices geometries.}
\label{fig:xcomp}
\end{figure}

\NEW{While the effect of gravity due to channel orientation helps reducing the pressure losses across the system, 
the short channels close to the flow inlet act as short circuit path,
allowing high mass flow rates of vapor directly from inlet to outlet. This has the clear disadvantage that high flow rates
of vapor cannot condense efficiently over a short distance. The described mass flow rate distribution has a strong effect
on the local vapor quality, as depicted in 
Fig.~\ref{fig:xcomp}. Generally, for a channel of fixed length, high flow rates correspond to
a high vapor quality at the discharge. This phenomenon is particularly evident in case "c", where the lower sub-channel
with the higher flow rate does not provide a good condensation due to its short length. The designer should seek for a
balanced distribution of the vapor qualities at the discharge of each channel in order to exploit best the heat transfer
area.}

\begin{figure}[h!]
\centering
\subfigure[device A]{\includegraphics[width=.7\linewidth]{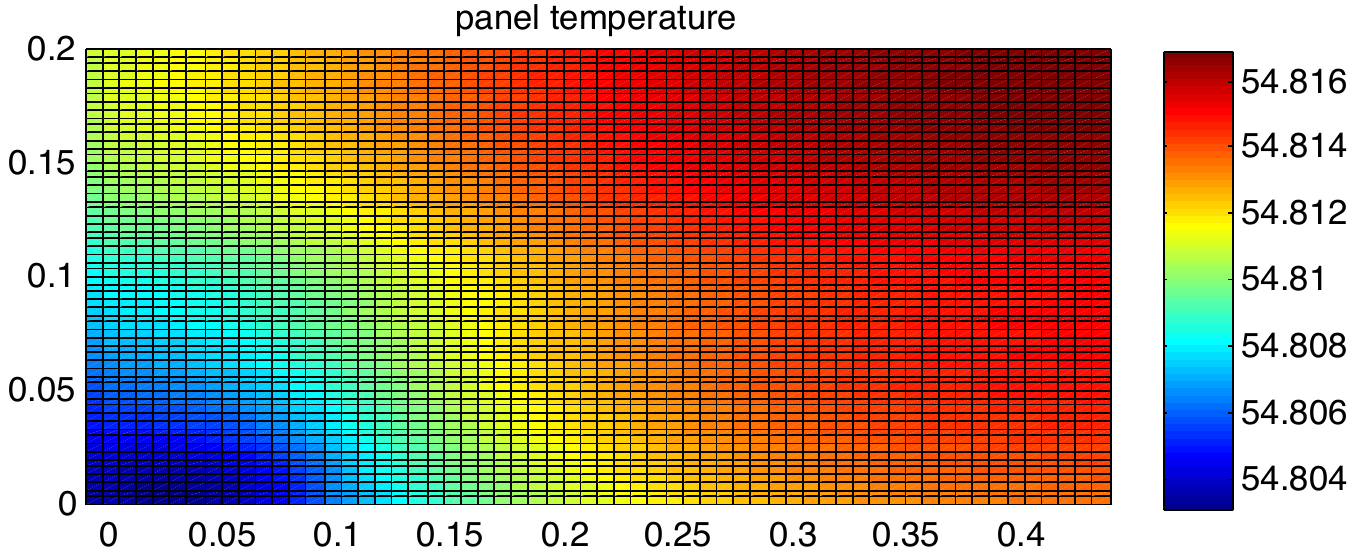}}
\subfigure[device B]{\includegraphics[width=.7\linewidth]{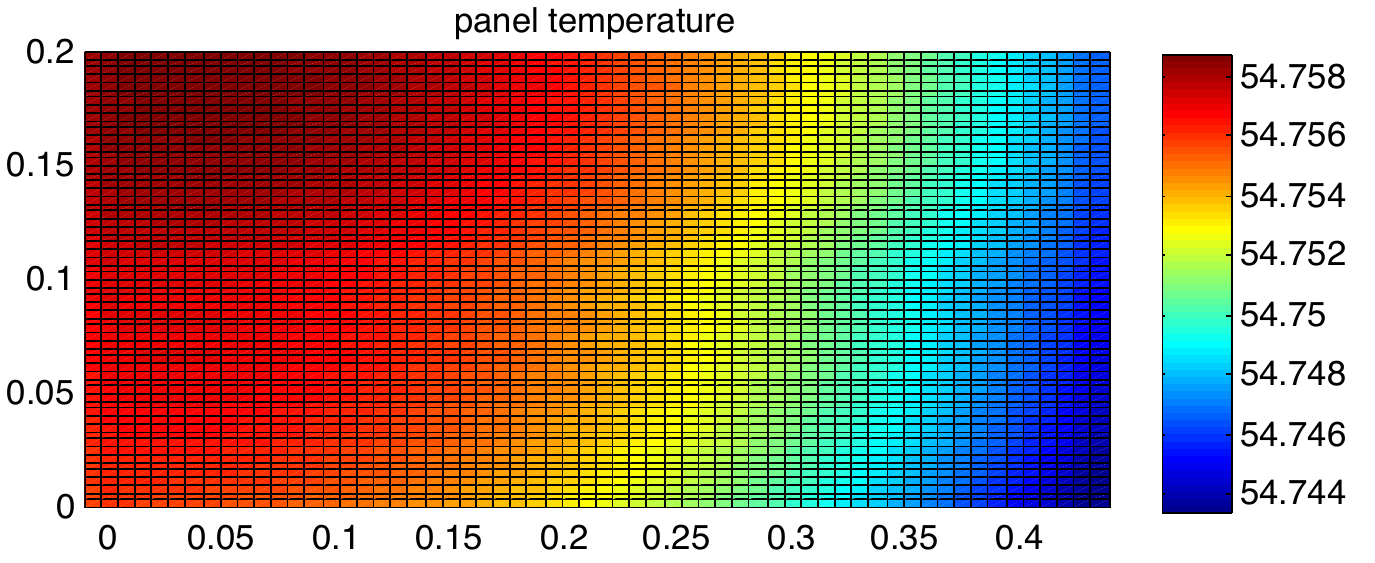}}
\subfigure[device C]{\includegraphics[width=.7\linewidth]{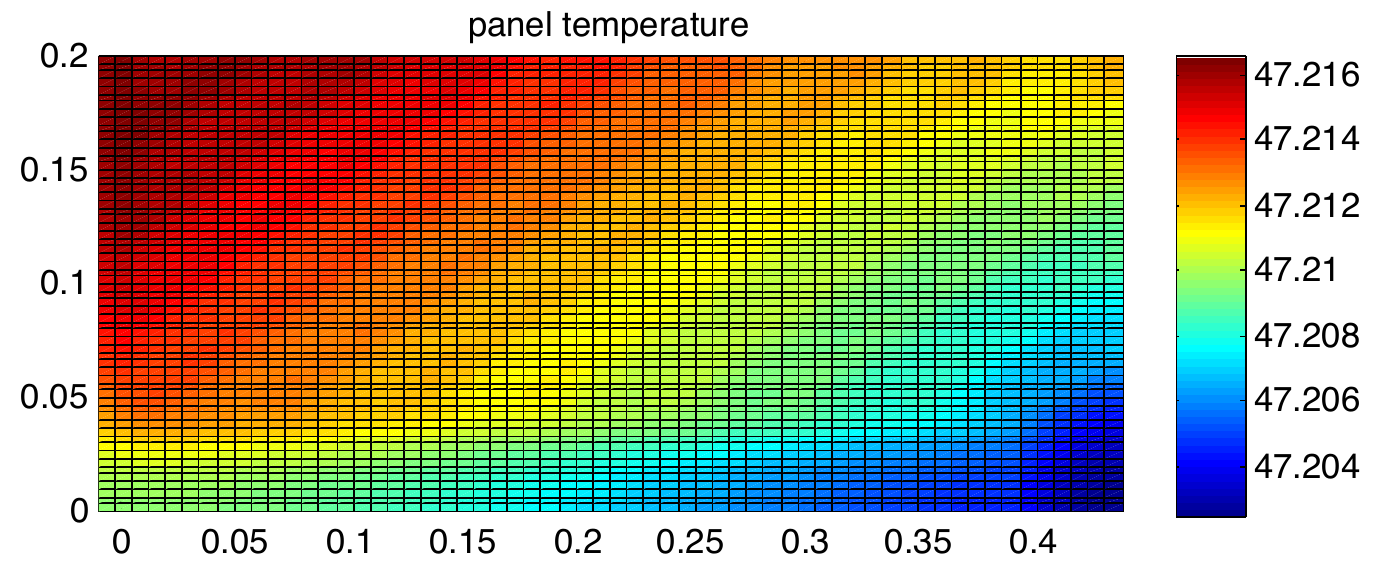}}
\caption{Comparison of panel temperature for the three device geometries.}
\label{fig:tcomp}
\end{figure}

\begin{figure}[h!]
\centering
\subfigure[device A]{\includegraphics[width=.7\linewidth]{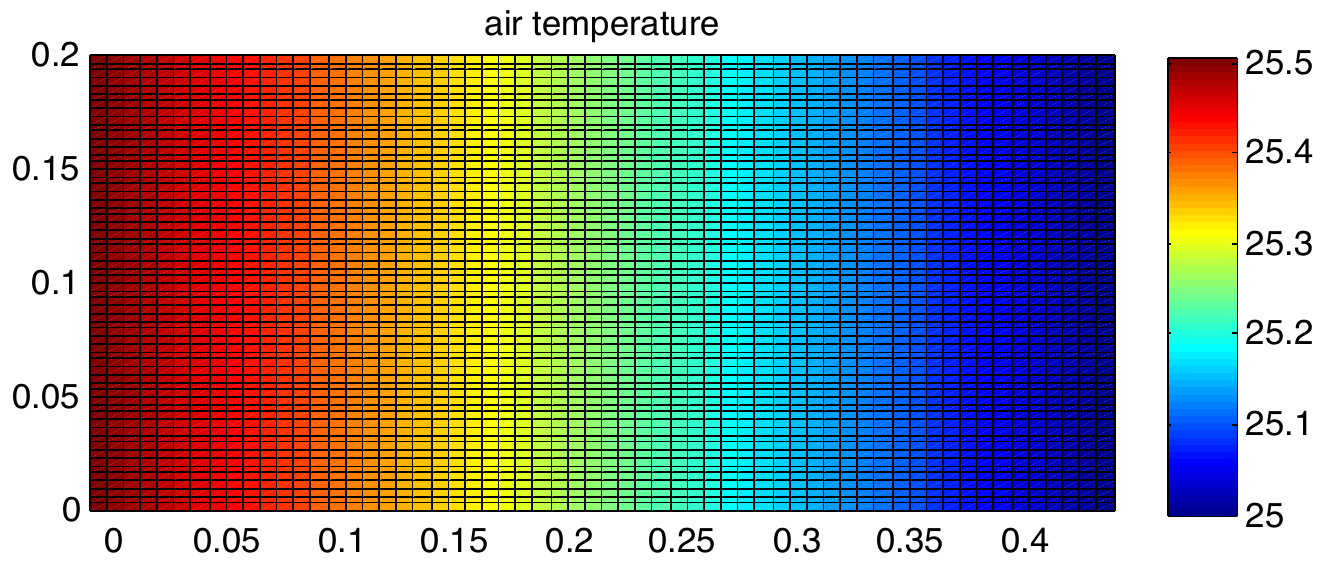}}
\subfigure[device B]{\includegraphics[width=.7\linewidth]{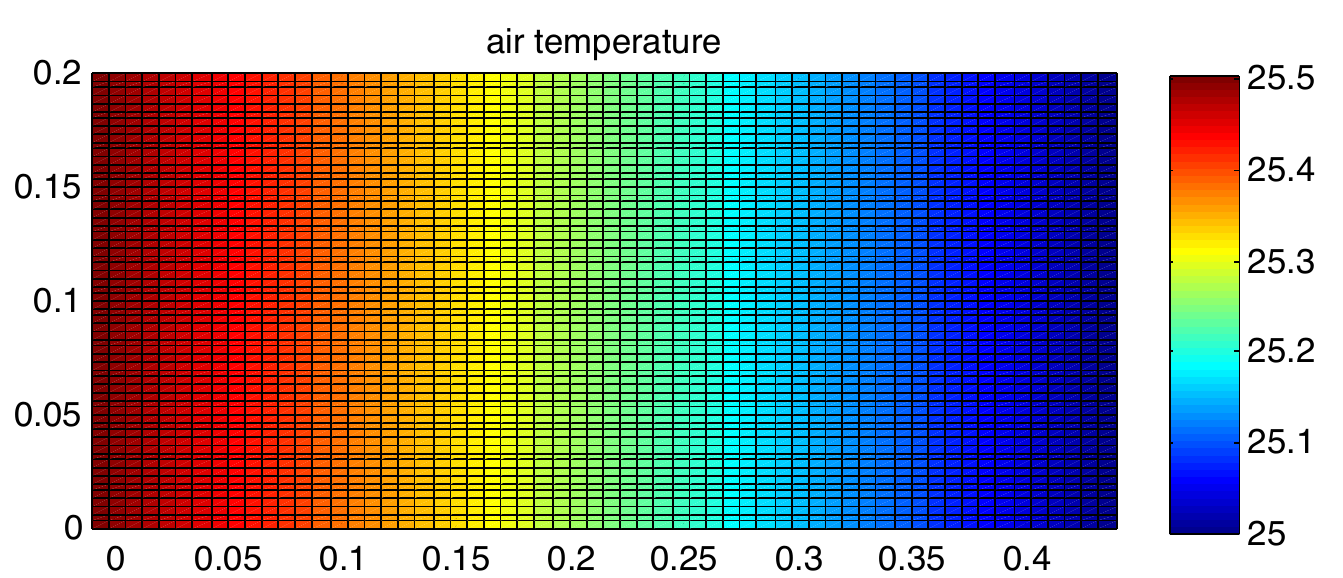}}
\subfigure[device C]{\includegraphics[width=.7\linewidth]{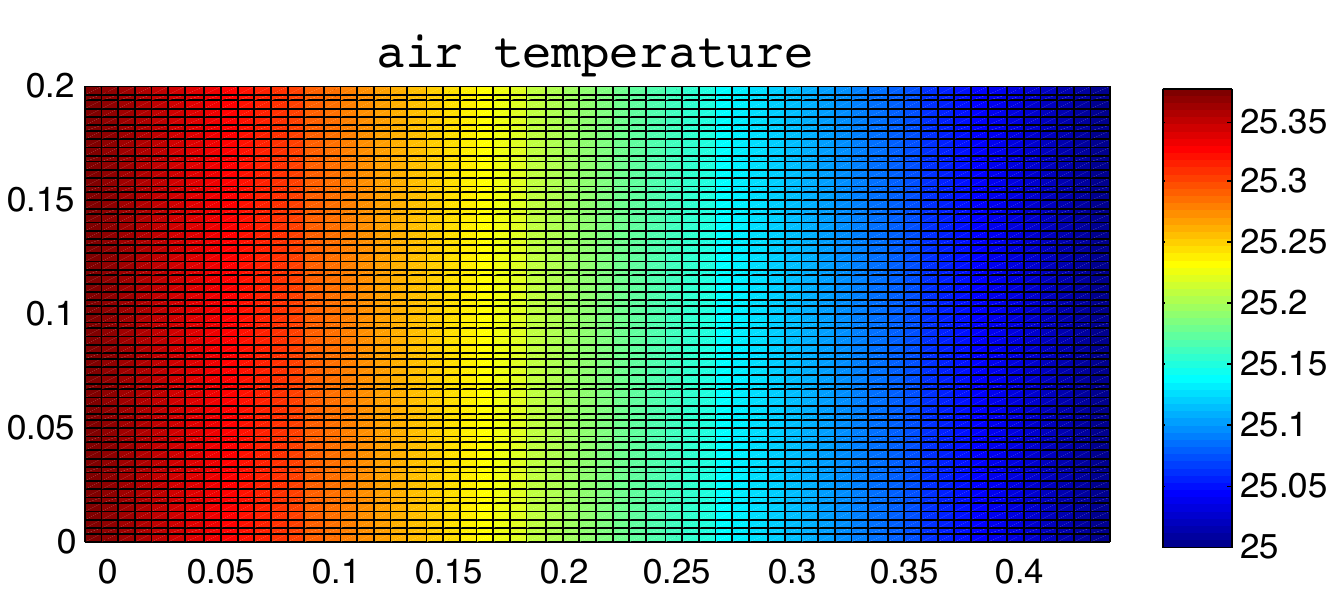}}
\caption{Comparison of air temperature for the three device geometries.}
\label{fig:tacomp}
\end{figure}

\NEW{Fig.~\ref{fig:tcomp} shows the value of the panel temperature for the three different geometries. 
Results indicate an almost constant temperature distribution. This is characteristic
of a two-phase system where the condensation heat transfer coefficients are orders of magnitude higher than those of
the air side. Fig.~\ref{fig:tacomp} shows the evolution of the air temperature for the three different geometries. This plot is a good
representation of the total heat transferred by the panel, representing the sensible heating of the air stream. 
Considering the original boundary condition of a fixed inlet mass flow rate of vapor, a higher air temperature difference 
indicates a higher amount of transferred heat. While case "a" and case "b" 
are comparable, case "c" shows a lower air temperature at the
discharge of the panel, clearly indicating a lower heat transfer to the air. This is well supported by the mass flow rate
distribution and vapor quality plots. 
We also can notice that in Fig.~\ref{fig:tacomp} 
the air temperature differences are smaller in case "c"
compared to "a" and "b".
This is probably to be ascribed to the 
fact that a mass flux is enforced as boundary condition in the simulation model (mass flow rate of vapor per unit area) and not power.}

\begin{figure}[h!]
\centering
\subfigure[device A]{\includegraphics[width=.7\linewidth]{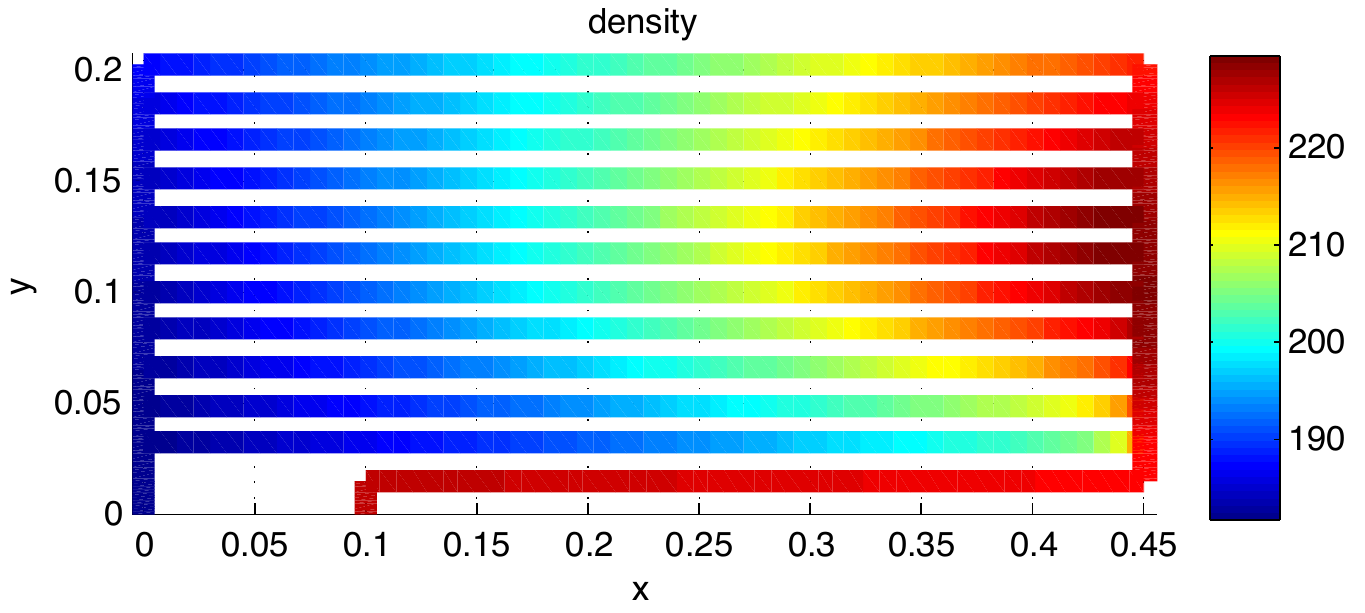}}
\subfigure[device B]{\includegraphics[width=.7\linewidth]{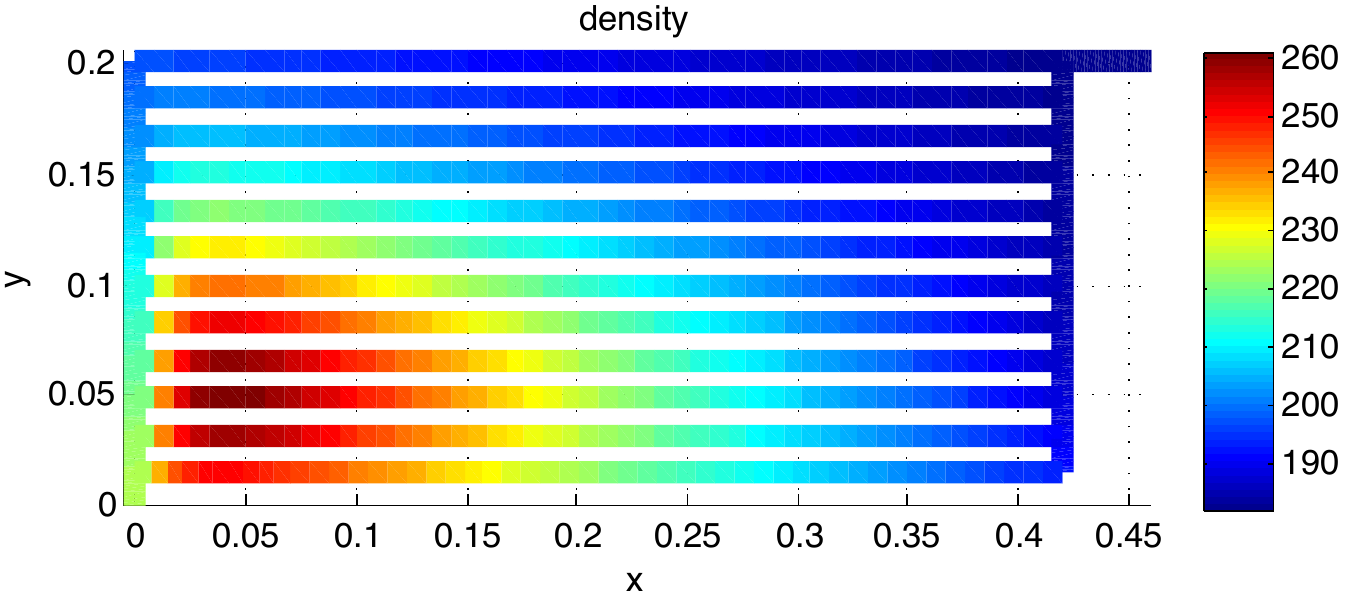}}
\subfigure[device C]{\includegraphics[width=.7\linewidth]{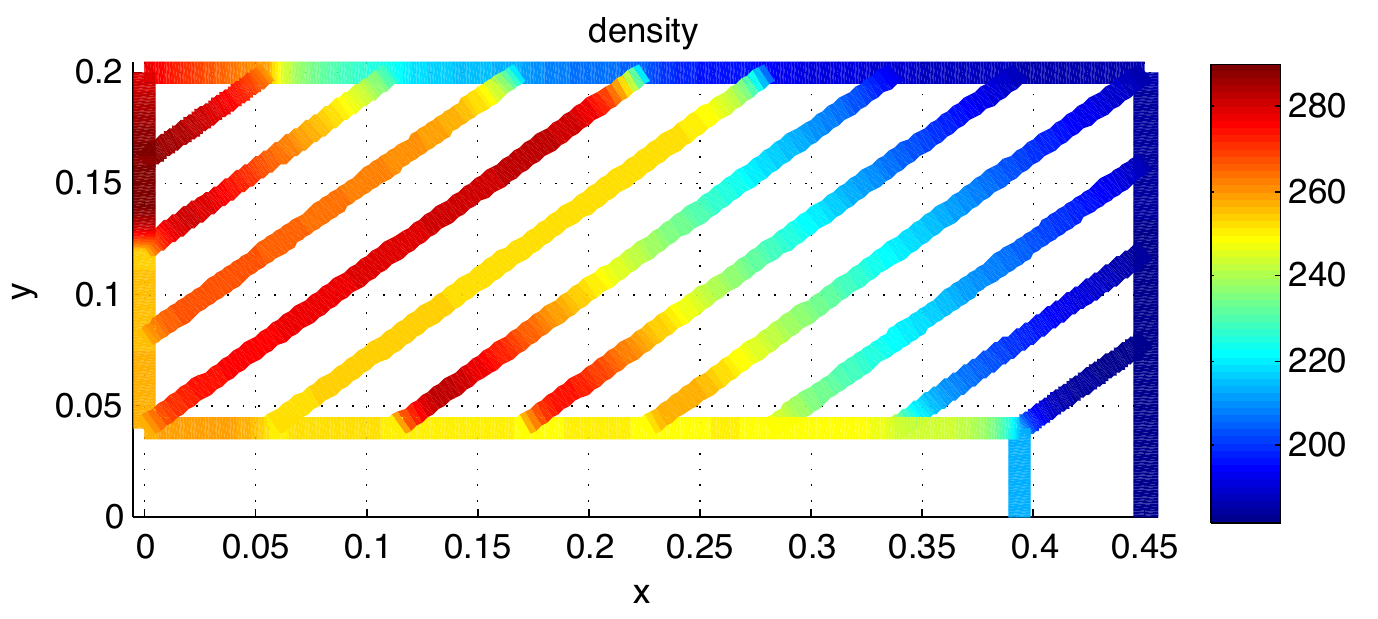}}
\caption{Comparison of the density of the two phase fluid for the 
three device geometries.}
\label{fig:dcomp}
\end{figure}

\NEW{Fig.~\ref{fig:dcomp} shows the spatial distribution of the so called 
two phase density or bulk density for the three different device geometries.
This quantity represents a weighted density between liquid and vapor 
densities, the weighting factor being the vapor quality. 
This means that the bulk density is the sum of the vapor and liquid densities multiplied by the vapor quality (vapor phase fraction) and its complement (liquid phase fraction), respectively. As a result, portions of the channels with higher bulk densities represent 
a fluid in a state with a higher content of liquid phase. 
To interpret Fig.~\ref{fig:dcomp} we can directly refer to 
Fig.~\ref{fig:xcomp}, so that high flow rates imply a high vapor quality at the discharge and relatively low two phase densities. As discussed 
for Fig.~\ref{fig:xcomp}, this
latter phenomenon is particularly evident in case "c" where the lower sub-channel with the higher flow rate does not provide a good condensation due to its short length and low vapor densities occur. As for the vapor quality, the design should seek for a balanced distribution of the two phase densities at the discharge of each channel in order to exploit at the best the heat transfer area and in order to have a balanced distribution of the liquid and vapor phases across the condensing panel.}

\subsection{Comparison with measured data}\label{sec:numres2}

In this section we carry out a set of simulation runs to validate the performance 
of the computational model \NEW{on realistic geometries and fluid-dynamical data.}
The experimental campaign and test set-up used for the validation follows closely what is presented in~\cite{experimental} and~\cite{iecon11}. As described in~\cite{iecon11}, the investigated cooling system is a thermosyphon device constituted of: an evaporator body, a vapor riser, a condenser (stack of roll-bonded panels) and a liquid downcomer. The evaporator can accommodate two ABB HighpakTM power semiconductor modules. Once the modules are in operation the evaporator collects the heat transferred by means of an evaporating fluid. The evaporator is designed \NEW{in such a way that} at its discharge the liquid is separated from the vapor. The liquid is brought back to the evaporator inlet while the vapor travels toward the condenser through the vapor riser. At the inlet of the condenser a vapor distributor feeds the stack of aluminum panels, equally distributing the mass flow among them. The panels are so-called roll-bonded panels, constituted of two aluminum sheets bounded together over almost the entire surface. Where this bounding is not present, a channel is generated, allowing the passage of the two-phase flow. The heat is rejected to the ambient by means of natural convection, the vapor is brought back to liquid conditions. Finally, the liquid \NEW{is} driven back to the evaporator inlet by gravity.
The same aluminum panels and stack geometrical layout as presented in~\cite{iecon11} is the subject of the investigation. The condenser is a stack of 13 panels 500 mm wide and 250 mm high, 1.2 mm thick, and equally spaced with a pitch of 18mm. Each panel contains 11 horizontal flow channels of a nominal length of 390 mm. The flow channel is formed on both sides of the panel with \NEW{isoscele trapezoidal} sections, the base and the height measuring 10 mm and 2.1 mm, respectively. The vapor and liquid phases are distributed to and collected from the panels by means of collectors of
19 and 16 mm internal diameter, respectively. Detailed drawings are available in~\cite{iecon11}, while a detailed description of the experimental measurement techniques is presented in~\cite{experimental}.
\NEW{The experimental conditions are summarized in Tab.~\ref{tab:tabdata}.}
\begin{table}[h!]
\centering
\begin{tabular}{|l | c|}\hline Fluid & R245fa \\ \hline
Refrigerant charge & 2 Kg \\ \hline
Filling ratio & 0.5 \\ \hline
Heat Load & 200 - 1600 W \\ \hline
Ambient temperature & 298.15 K \\ \hline
Air cooling regime & Natural convection \\ \hline
\end{tabular}
\caption{Experimental conditions.}
\label{tab:tabdata}
\end{table}

\begin{figure}[h!]
\begin{center}
\includegraphics[width=0.95\columnwidth]{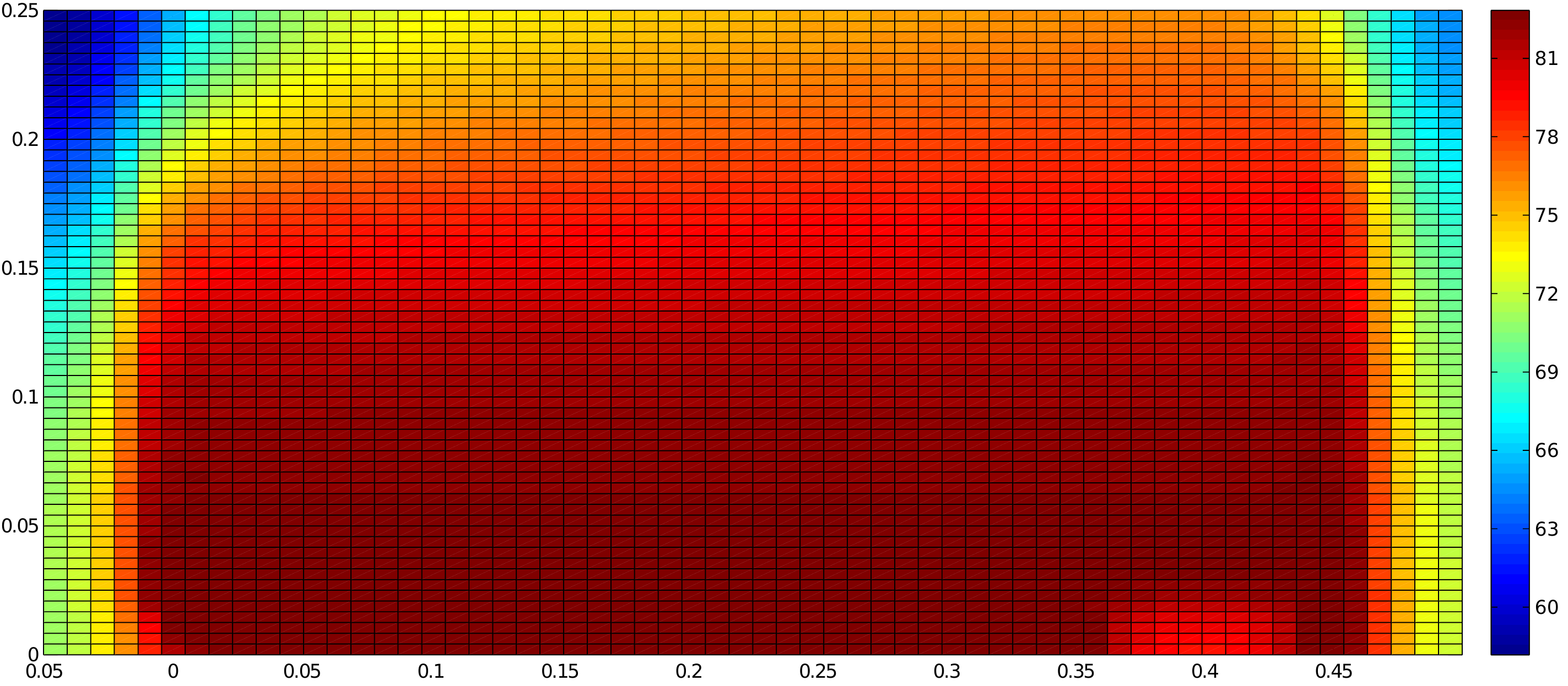}
\end{center}
\caption{Computed panel temperature for a total dissipated power of
  $1500$\,W.} 
\label{fig:panel_num}
\end{figure}

Figure~\ref{fig:panel_num} presents the computed panel temperature corresponding to a power inflow of 1500 W, air inlet temperature of 25 $^{\circ}$C and natural convection operation. It is observed that the panel is almost isothermal. This is characteristic of the investigated system. A condensing fluid in the panel channels is characterized by high heat transfer coefficients, orders of magnitude higher than heat transfer coefficients typical of natural convection in air. The heat transfer conditions as well as the nature of the panel, sufficient thickness, small distance between channels and relatively high thermal conductivity of the aluminum result in an almost constant panel temperature.
An almost constant temperature of the condenser panel is what we are looking for from an application point of view. It allows to overcome a common drawback of a standard heat-sink based system, where the metallic fin does not behave as a perfect fin (constant temperature) but has a temperature gradient from base to tip, resulting in a limited efficiency. Having an almost constant temperature results in an efficiency of the fin close to unity.
The panel border is the coldest part. The low temperature in this region is is due to boundary effects. While the rest of the panel has an almost constant temperature, we can still identify a hotter region in the lower part of the panel compared to the top part.
\begin{figure}[h!]
\begin{center}
\includegraphics[width=0.95\columnwidth]{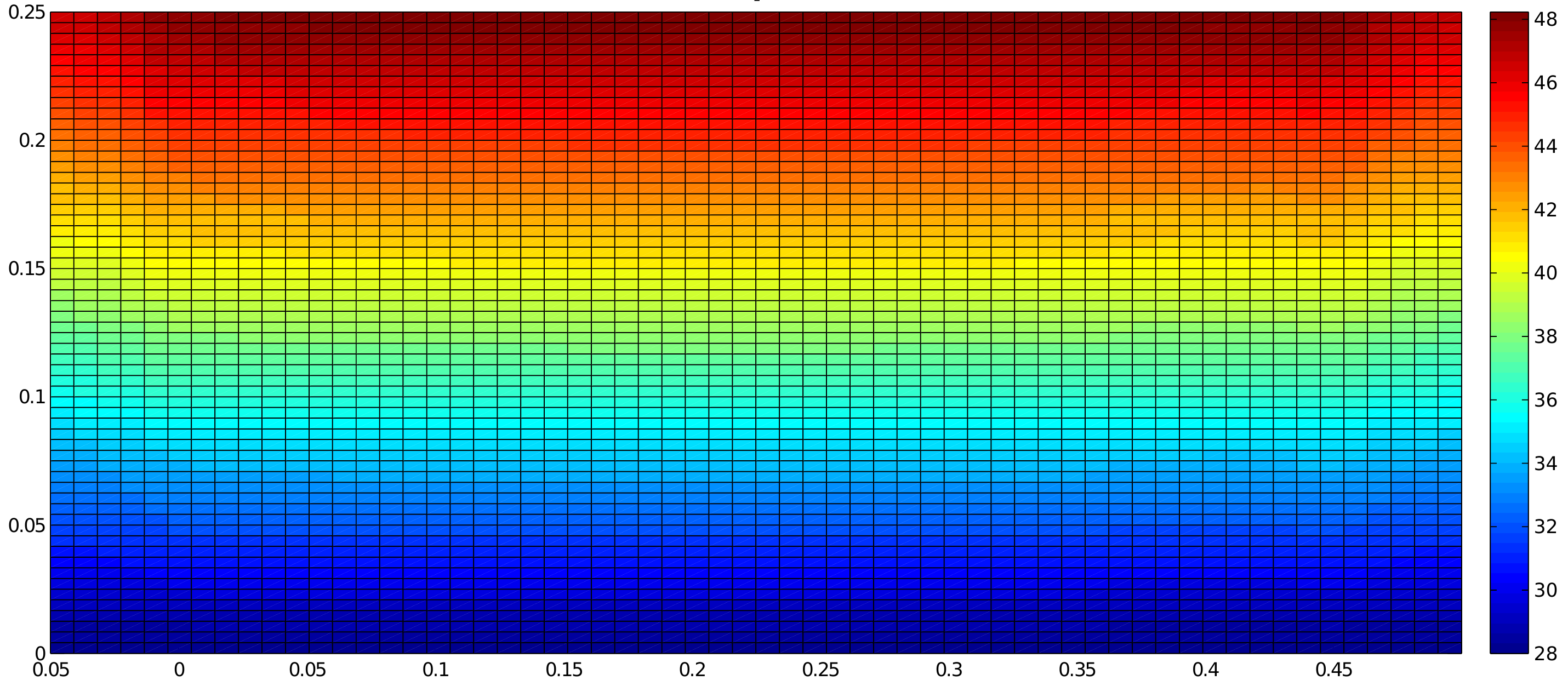}
\end{center}
\caption{Computed air temperature for a total dissipated power 
of $1500$\,W.}
\label{fig:air_num}
\end{figure}

Figure~\ref{fig:air_num} depicts the air temperature between two panels. The values are averaged in the direction perpendicular to the panel surface. The temperature pattern is characteristic of the transfer of sensible heat from panel to air in natural-convection operation. The \NEW{large temperature difference} between inlet and outlet of the condenser panel \NEW{results} from the \NEW{small air velocity} typical of natural convection. The maximum allowed temperature difference between inlet and outlet is usually a design parameter, and the designer of the device tries to optimize the system in order to match this value. 
\NEW{A higher allowed temperature difference makes} 
it possible to shrink the size of the device. On the other hand, when required, a decrease of the maximum temperature difference can be obtained by increasing the number of panels or the panel area. Since the panel temperature decreases from bottom to top, while the air temperature increases in that direction, the temperature difference between air and panel is largest at the panel bottom. This means that the heat flux from panel to air is maximum at the panel bottom.

The mass flow rate distribution among the panel channels plays an important role in the behavior of the condenser. During operation, the two-phase 
fluid \NEW{tends to} flow in the horizontal channels \NEW{suitably} paralleled. Considering the fact that the flow path through the panel and hence the flow resistance is smallest for the bottom channel, a decrease in mass flow rate from bottom to top is expected. A certain inhomogeneity in mass flow rate must therefore always be accounted for in the type of parallel connection of channels. Due to the higher mass flow rate, and consequently higher velocity in the bottom channel, a lower fluid residence time per channel length results. Consequently, it is expected that a longer channel length is needed to complete condensation.
\begin{figure}[h!]
\begin{center}
\includegraphics[width=0.95\columnwidth]{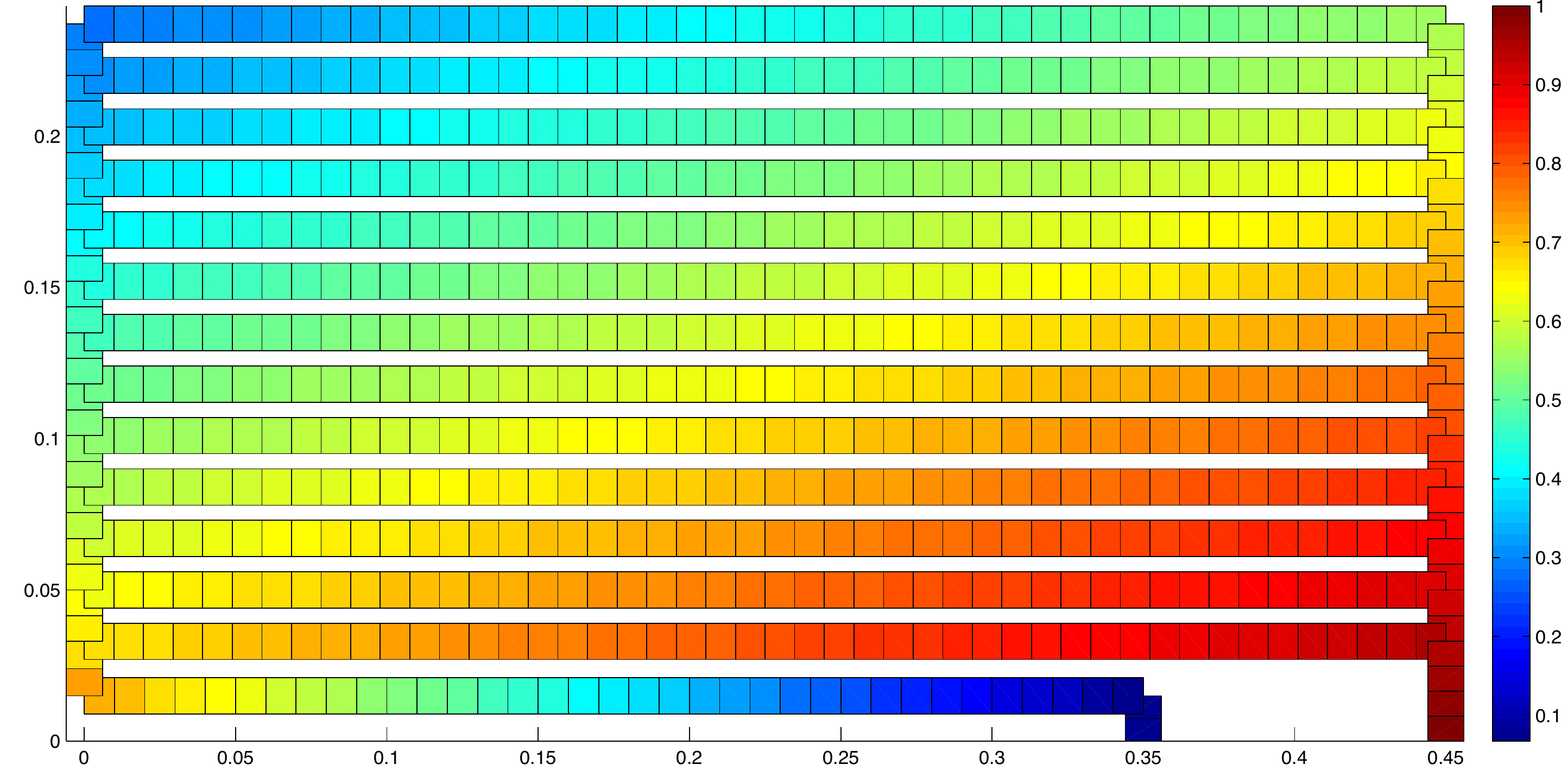}
\end{center}
\caption{Computed refrigerant vapor quality for a total dissipated
  power of $1500$\,W.} 
\label{fig:x_num}
\end{figure}

This is indeed observed in the simulation results in Figure~\ref{fig:x_num}, 
showing the local vapor quality in the channels. For the bottom channel, a longer distance from the channel inlet is needed for the vapor quality to decay to a certain value. Consequently, the vapor quality at the channel end, i.e. at the left in the figure, is highest for the bottom channel and lowest for the top channel. Furthermore, from the energy balance, it is clear the condensation of the highest mass flow rate in the bottom channel requires the largest heat flow rate from channel to air. Since all channels have the same surface area, one expects the heat flux to be highest for the bottom channel and lowest for the top channel. This \NEW{closely agrees with the} observation of maximum temperature difference between panel and air at the bottom and the 
\NEW{corresponding} maximum heat flux between panel and air in this region.

The designer may try to minimize the observed differences in performance between
the condenser channels by optimizing the channel design. For example, one may
try to achieve the same vapor quality at the end of all condenser channels.
Complete condensation and hence low vapor quality is fundamental to guarantee a
safe and reliable operation of the device, since a re-wetting of the evaporator
surface is mandatory. It is exactly this kind of optimization tasks for which the present mathematical model is beneficial, as it provides insight in the detailed performance and behavior of the cooling device.
\begin{figure}[h!]
\begin{center}
\includegraphics[width=0.99\columnwidth]{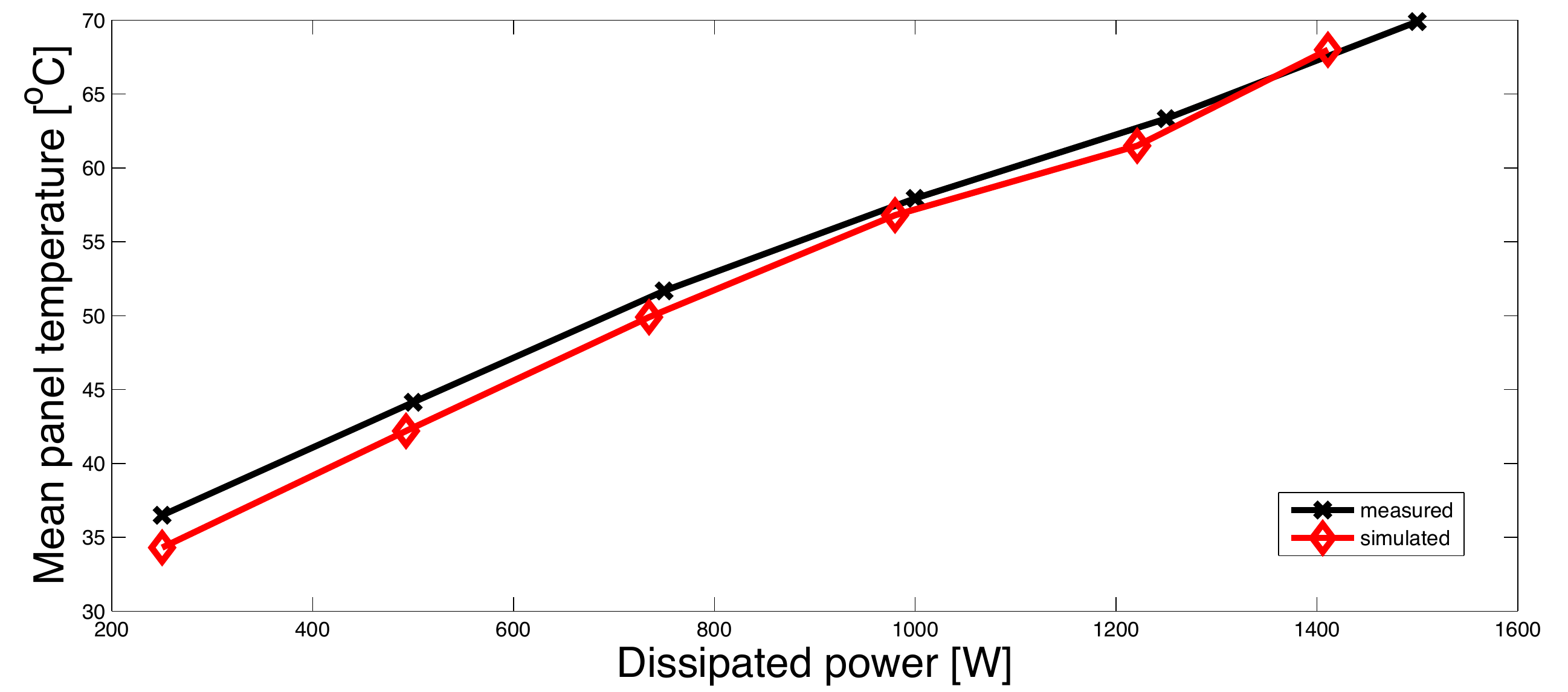}
\end{center}
\caption{Average panel temperature (computed vs. measured) as a
  function of the total dissipated power.}
\label{fig:mean_temp}
\end{figure}

Figure~\ref{fig:mean_temp} shows a plot of the mean temperature of the panel as a function of the dissipated power. While the 
\NEW{computed temperature distribution describes in great detail the operation of the device, the mean panel temperature is a 
synthetic parameter} for the designer to validate \NEW{in an immediate manner} 
the predictive capabilities of the code. \NEW{Agreement of numerical results of Figure~\ref{fig:mean_temp} with experimental data 
is striking and indicates that, although based on many simplifying assumptions, our model does have 
very good predictive accuracy.}

\section{Conclusions and Future Work}\label{sec:conclusions}

In this article we have proposed and numerically implemented
a multiscale thermo-fluid mathematical model for the description of a 
condenser component of a novel two-phase thermosyphon cooling 
system presented in~\cite{experimental,iecon11}. 
The condenser consists of a set of roll-bonded vertically mounted fins among which air flows by either natural or forced convection and plays 
an important role in the industrial design of advanced power 
electronics systems.

The mathematical model \NEW{is} developed with the aim of
deepening the understanding of the various thermo-fluid 
mechanisms that determine the performance of the condenser 
in view of a further optimization of the cooling device.
The adopted approach is based on a multiscale formulation 
meant to reduce as much as possible the complexity required 
by a fully three-dimensional (3D) simulation code while 
maintaining reasonable predictive accuracy.

More specifically, the flow of the two-phase coolant within the 
condenser fins is modeled as a 1D network of pipes, while heat 
diffusion in the fins and its convective transport in the air 
slab are modeled as 2D processes.
The resulting mathematical problem consists of a system of nonlinearly
coupled PDEs in conservation form that are characterized by a mixed parabolic-hyperbolic character with possible presence of strongly advective dominating terms. A fixed point iterative map is used 
to reduce the computational effort to the successive solution of
a sequence of decoupled linear stationary boundary value problems in
the 1D channel pipe network and in the 2D air domain, respectively.

For the numerical approximation of the above differential problems
a Primal Mixed Finite Element discretization method with upwind 
stabilization is used for the 1D coolant flow while a 
Dual Mixed-Finite Volume scheme with Exponential Fitting 
stabilization is used for 2D heat diffusion and convection. 

Extensive numerical tests are carried out to validate the stability
and accuracy of the proposed schemes on several benchmark problems
whose solution is characterized by the presence of steep interior and 
boundary layers. The obtained results confirm the good accuracy
of the proposed formulation and its \NEW{ability in} satisfying 
a discrete maximum principle. This latter property confers robustness
to the simulation tool and makes it suitable for heavy duty
use in industrial applications.

The solver is then thouroughly applied to the numerical study
and parametric characterization of a two-phase coolant system
with realistic industrial geometry. The output of the simulations
provide a complete map of the principal thermal and fluid dynamical 
variables of the problem (air temperature, coolant fluid pressure 
and vapor quality) that are extensively used by the project
engineer to quantitavely design a novel device structure.
Two groups of simulations are performed for \NEW{the validation of the computational
algorithm.}
In a first set of runs, the code is used to analyze the impact of channel geometry on the distribution of mass flow rate, vapor quality and panel temperature.
In a second set of runs, the simulated average panel temperature of a given
realistic cooler geometry is compared with available experimental data.
Despite the several simplifying model assumptions introduced in 
the condenser mathematical description, the obtained results turn out 
to be in very good agreement with measures thus providing a sound 
indication of model reliability. 

\NEW{Even if} applied to a problem arising in a 
specific area of thermo-fluid dynamical industrial applications,
the multiscale modeling approach proposed in the present work
\NEW{can be used to study problems arising in other scientific contexts.
For example, the computational
model to couple 2D heat convection-diffusion and 1D flow in a pipeline network} shares a close resemblance
with the mathematical and numerical treatment of flow and mass 
transport in biological tissues that has been recently investigated
in~\cite{dAngelo2007,Shipley2010,Erbertseder2012} 
and references cited therein.
This interesting similarity might be profitably used 
\NEW{to apply to these latter novel bio-technological applications} 
solution methods that \NEW{in this article are proved to enjoy properties
of accuracy, stability and conservation.}

Further research activity will be devoted to the:
\begin{itemize}
\item topological optimization of the channels layout; 
\item integration of the condenser model in a complete thermosyphon loop simulation tool including evaporator body and connections;
\item analysis of the existence of a fixed point 
of the iterative map and its possible uniqueness.
\end{itemize}


\appendix
\section{Dimensionality Reduction 
of the Heat Con\-vec\-tion-Dif\-fu\-sion Equations}
\label{sec:3d_2_2d}

\begin{figure}[h!]
\begin{center}
\includegraphics[width=.8\linewidth]{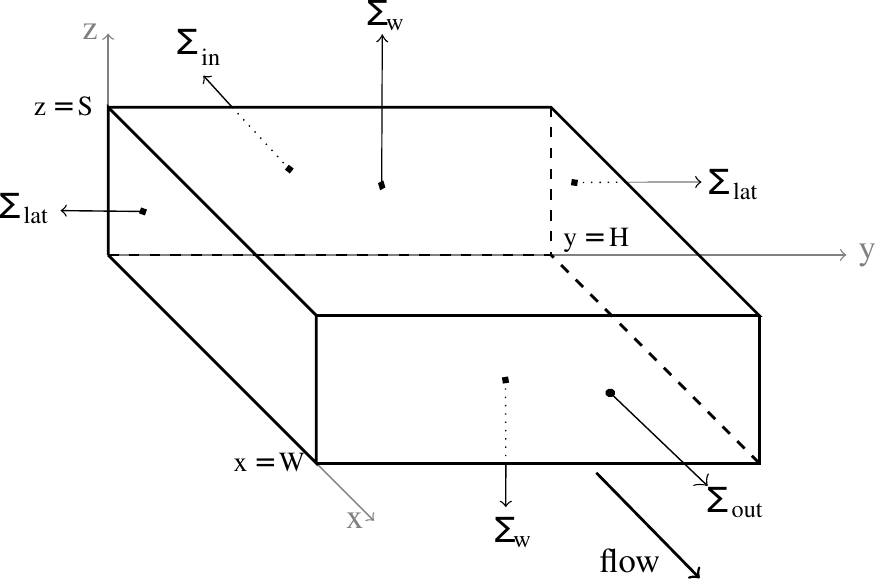}
\end{center}
\caption{Three-dimensional domain for the heat 
convection-diffusion model problem.}
\label{fig:air_domain}
\end{figure}

In this section we illustrate the model reduction procedure 
that allows to derive under the assumptions (H1)--(H7) of 
Sect.~\ref{sec:mathmod2d} the simplified 2D 
model~\eqref{eq:2d_air_panel} from the corresponding 3-dimensional
equations for heat convection and diffusion in the condenser walls 
and in the air between two plate walls. 
In order to describe the dimensionality reduction procedure, 
we start from
the following model problem set in the 3D computational domain 
$\Omega$ depicted in Fig.~\ref{fig:air_domain}:
\begin{equation}\label{eq:3dmodelproblem}
\begin{array}{rcll}
\nabla \cdot \left(-k \nabu u +\rho c \vector{v} u\right) &=& 0 & 
\qquad \mbox{in } \Omega \\[3MM]
u &=& u_{in} & \qquad \mbox{on } \Sigma_{in} \\[3MM]
-k \nabu u \cdot \vector{n} &=& 0 & \qquad \mbox{on } \Sigma_{out} \\[3MM]
\left(-k \nabu u + \rho c \vector{v} u\right) \cdot {\vector{n}} &=& 
\textit{h}\ (u - u_{w}) & \qquad \mbox{on } \Sigma_w \\[3MM]
\left(-k \nabu u + \rho c \vector{v} u\right) \cdot \vector{n} &=&  0 
& \qquad \mbox{on } \Sigma_{lat}
\end{array}
\end{equation}

We notice that the model problem~\eqref{eq:3dmodelproblem} 
may describe either the forced heat convection between two fins or
the heat diffusion in one fin wall. 
Referring to Fig.~\ref{fig:air_domain}, 
in the latter case, we have
$\Sigma_{in} = \Sigma_{out}=\emptyset$, $\vector{v}=\vector{0}$ 
and $\Sigma_{lat}=\displaystyle \cup_{i=1}^4 \Sigma_i$, while
in the former case we have 
$\Sigma_{in} = \Sigma_1$, $\Sigma_{out}=\Sigma_3$ and 
$\Sigma_{lat}=\Sigma_2 \cup \Sigma_4$. 
The {\it contact walls} 
$\Sigma_w$, located at $z=0$ and $z=S$ respectively, 
represent the boundaries where heat exchange occurs.
According to assumption (H4), the {\it convection velocity} 
$\vector{v}$ is directed along 
the $x$ axis, so that it can be expressed as
\begin{equation}
\label{eq:velocityansatz}
\vector{v}(z) = \vector{V}\ B(z) 
\end{equation}
where $\vector{V}$ is a constant vector directed along the $x$ axis and $B(z)$ is a dimensionless scalar shape function accounting for the velocity boundary layer in the $z$ direction.

The unknown function $u = u(x,y,z)$ represents a temperature 
(either air temperature or wall temperature), $\rho$ is the density of 
the medium contained in the domain $\Omega$, $c$ is the specific heat 
capacity of the medium and $k$ is its thermal conductivity. 

Temperature is fixed at the inlet surface $\Sigma_{in}$ to a 
given value $u_{in}$.
On the contact surfaces $\Sigma_w$ the outflow heat flux is 
proportional to the difference between temperature $u$ and 
the wall temperature $u_w$, through the {\it heat transfer coefficient} 
$\textit{h}$. $\vector{n}$ is the outward unit vector along the external 
surface of the domain. 
 
According to assumptions (H1) and (H3), 
the conditions at the upper and lower contact surface $\Sigma_w$ 
are symmetric. Therefore, we can define an adiabatic plane at 
$z=S/2$ which allows us to consider only the portion of space between 
the adiabatic surface and one of the two contacts $\Sigma_w$, 
for example that located at $z=0$.

We start our dimensionality reduction procedure
by assuming the following ansatz for the unknown $u$
\begin{equation}
\label{eq:sep_var}
u(x,y,z) = U(x,y) \, Z(z),
\end{equation}
where $U=U(x,y)$ expresses temperature variation 
in the $xy$ plane, while $Z=Z(z)$ is a dimensionless shape 
function accounting for temperature variation between the contact 
surface and the adiabatic plane located at $z=S/2$. 
The separated variable form of temperature 
distribution~\eqref{eq:sep_var} agrees well with assumptions 
(H6) and (H7) according to which a mild variation of temperature
between two neighbouring contact surfaces is to be 
expected. The next step consists in examining the dependence 
of problem coefficients on the unknown $u$.
The heat capacity $c$ can be taken as a constant~\cite{wassermann1966}.
The same holds for the density $\rho$ (cf. assumption (H5)).
As far as the thermal conductivity $k$, the following power 
law can be used~\cite{wassermann1966}
\begin{equation}\label{eq:power_law_kappa}
k(x,y,z) = k_0 \left( \Frac{u(x,y,z)}{u_0} \right)^\beta =
k_0 \left( \Frac{U(x,y) Z(z)}{u_0} \right)^\beta
\end{equation}
where $k_0$, $u_0$ and $\beta$ are suitable constants.

Integration of the balance equation 
in the vertical direction and the use 
of~\eqref{eq:velocityansatz},~\eqref{eq:sep_var} and~\eqref{eq:power_law_kappa} yield
\begin{equation}\label{eq:NS_sp}
\lambda_1 \grad_{xy} \cdot \left( -k_0 \left( 
\Frac{U}{u_0}\right)^\beta \grad_{xy} U \right)
+ \lambda_2 \grad_{xy} \cdot \left(\rho c \vector{V} U \right)
+ \mathbb{I} = 0
\end{equation}
where 
\begin{equation}\label{eq:lambda}
\lambda_1 := \dint_{0}^{S/2} Z^{\beta+1} (z)\ dz, \qquad
\lambda_2 := \dint_{0}^{S/2} Z (z)\ B(z)\ dz
\end{equation}
and 
$$
\mathbb{I}: = \dint_{0}^{S/2} \partial_z 
\left( -k \partial_z u \right) dz,
$$
while $\grad_{xy} (\cdot)$ is the gradient operator 
with respect to the directions $x$ and $y$ only, 
The quantity $\lambda_1$ modulates the
variation of thermal conductivity in the $z$ direction while the
quantity $\lambda_2$ is related to the 
shape of the thermal boundary layer arising at the 
interface between air and panel. 
Using Gauss theorem to treat the quantity $\mathbb{I}$ we get
$$
\mathbb{I}= \left[ - k \partial_z u \right]_{0}^{S/2} = 
- \left(- k \left. \partial_z u \right|_{z=0} \right)
$$
because $-k \left. \partial_z u \right|_{z=S/2} = 0$ 
under the assumption of adiabatic surface, so that
we can rewrite condition~\eqref{eq:3dmodelproblem}$_4$ as
$$
\textit{h} \left( u|_{z=0} - u_w \right) = 
(-k \nabu u + \rho c u \vector{v}) \cdot \vector{n} \Big|_{\Sigma_w} = 
-k \nabu u \cdot \vector{n} \Big|_{\Sigma_w} = 
- \left(- k \partial_z u|_{z=0} \right) = \mathbb{I}.
$$
Therefore, upon rescaling the shape function $Z$ in such a way that 
$Z(0) = 1$, equation~(\ref{eq:NS_sp}) becomes
$$
\lambda_1 \grad_{xy} \cdot \left( - k_{xy} \grad_{xy} U \right) 
+ \lambda_2 \grad_{xy} \cdot \left(\rho c \vector{V} U \right)
+ \textit{h} \left(U - u_w \right) = 0
$$
where we have defined the heat conductivity in the $xy$ plane
(contact surface)
$$
k_{xy}:= k_0 \left( \Frac{U(x,y)}{u_0} \right)^\beta
= k_0 \left( \Frac{U(x,y) Z(0)}{u_0} \right)^\beta = k(x,y,0).
$$
To end up with a 2D reduced model for heat convection and diffusion,
we need specify the exponent $\beta$. At low pressures
we typically have $\beta = 0.9$~\cite{wassermann1966}
so that $k$ is approximately a linear function of temperature.
This latter quantity
is used as a fitting parameter in the numerical simulations
reported in Sect.~\ref{sec:results}.
Thus, omitting the subscript $xy$ in the notation, and writing 
$u$ instead of $U$, the reduced 2D version 
of~\eqref{eq:3dmodelproblem} reads:
\begin{equation}
\label{eq:2d_air_reduced}
\left\{ \begin{array}{rcll}
\nabu \cdot (-k  \nabu u  + \rho {c} 
\widehat{\vector{v}} u) 
+ \widehat{\textit{h}} (u - u_w)  &=& 0 & \quad
(x,y) \in \Omega,\\
u &=& u_{in} &  \quad x=0, \\
-k  \nabu u \cdot \vector{n}  &=& 0 & \quad x = W,\\
(-k  \nabu u  + \rho c \widehat{\vector{v}} u) \cdot \vector{n} &=& 0 & 
\quad y=0, \ y=H,
\end{array}
\right.
\end{equation}
where $\Omega:=(0,W)\times(0,H)$, 
$\widehat{\textit{h}}:= \textit{h}/\lambda_{1}$, 
$\widehat{\lambda}:= \lambda_{2}/\lambda_{1}$
and $\widehat{\vector{v}}:= \widehat{\lambda} \vector{V}$.
Notice that the two heat balance equations~\eqref{eq:2d_air}
and~\eqref{eq:2d_panel} 
are special instances of~\eqref{eq:2d_air_reduced} 
upon setting $u=T_a$, $k=k_a$, $\rho=\rho_a$, $c=c_p$, 
$\widehat{\vector{v}}=\widetilde{\vector{v}}_a$, 
$\widehat{\textit{h}}= \widetilde{\textit{h}}_{aw}$ and $u_w = T_w$
in the case of the air temperature model, and
$u=T_w$, $k=k_w$, $\widehat{\vector{v}}=\vector{0}$, 
$\widehat{\textit{h}}={\textit{h}}_{aw}^\ast 
+ \textit{h}_{wc}^\ast$, and $u_w = 
(\textit{h}_{aw}^\ast T_a + \textit{h}_{wc}^\ast T_c)/
(\textit{h}_{aw}^\ast + \textit{h}_{wc}^\ast)$ 
in the case of the panel temperature model, respectively.

\section*{Acknowledgements}

RS was supported by the M.U.R.S.T. grant nr. 
200834WK7H005 {\em Adattivit\`a Numerica e di Modello
per Problemi alle Derivate Parziali}. 
CdF's work was partially funded by the {\em Start-up Packages and PhD Program project}, co-funded by Regione Lombardia through the {\em Fondo per lo sviluppo e la coesione 2007-2013}, formerly FAS program.

\bibliographystyle{plain}
\bibliography{paperabb}

\end{document}